%% file: main.tex
\documentclass[final,hidelinks,onefignum,onetabnum]{siamart220329}

\usepackage{amsfonts}
\usepackage{amsmath}
\usepackage{amssymb}
\usepackage{graphicx}
\usepackage{epstopdf}
\usepackage{algorithm}
\usepackage[italicComments=true, rightComments=false]{algpseudocodex}
\usepackage{caption}
\usepackage{subcaption}
\usepackage{booktabs}
\usepackage{multirow}
\usepackage{siunitx}
\usepackage{xspace}
\usepackage{mathtools}

\usepackage{color}
\ifpdf
  \DeclareGraphicsExtensions{.eps,.pdf,.png,.jpg}
\else
  \DeclareGraphicsExtensions{.eps}
\fi

\newcommand{\HPS}{Hierarchical Poincaré-Steklov\xspace}
\newcommand{\ignore}[1]{}
\newcommand{\pforest}{{\tt p4est}\xspace}

\renewcommand{\refeq}[1]{(\ref{#1})}  

\newcommand{\refsec}[1]{Section \ref{#1}}
\newcommand{\reffig}[1]{Figure \ref{#1}}
\newcommand{\refalg}[1]{Algorithm \ref{#1}}
\newcommand{\reftab}[1]{Table \ref{#1}}

\newcommand{\donna}[1]{{\color{blue}DC : #1}}
\newcommand{\amrex}{AMReX\xspace}
\newcommand{\eqn}[1]{(\ref{eq:#1})}
\newcommand{\Alg}[1]{Algorithm \ref{alg:#1}\xspace}

\newcommand{\real}{\ensuremath{\mathbb R}\xspace}
\newcommand{\DtN}{Dirichlet-to-Neumann\xspace}

\newcommand{\Ttau}{\ensuremath{\mathbf T^{\tau}}\xspace}
\newcommand{\Stau}{\ensuremath{\mathbf S^{\tau}}\xspace}

\newcommand{\Xtau}{\ensuremath{\mathbf X^{\tau}}\xspace}

\newcommand{\Gop}{\ensuremath{\mathbf G}\xspace}

\newcommand{\Lhtau}{\ensuremath{\mathcal L_\text{h}^{\tau}}\xspace}
\newcommand{\Atau}{\ensuremath{\mathcal A^{\tau}}\xspace}
\newcommand{\Linhtau}{\ensuremath{\mathcal L_\text{inh}^{\tau}}\xspace}

\newcommand{\Bctau}{\ensuremath{\mathcal B^\tau}\xspace}

\newcommand{\Talpha}{\ensuremath{\mathbf T^{\alpha}}\xspace}  
\newcommand{\Tbeta}{\ensuremath{\mathbf T^{\beta}}\xspace}  
\newcommand{\Tgamma}{\ensuremath{\mathbf T^{\gamma}}\xspace}  
\newcommand{\Tomega}{\ensuremath{\mathbf T^{\omega}}\xspace}  
\newcommand{\Ti}{\ensuremath{\mathbf T^{i}}\xspace}

\newcommand{\Talphap}[2]{\ensuremath{\Talpha_{#1#2}}\xspace}  
\newcommand{\Tbetap}[2] {\ensuremath{\Tbeta_{#1#2}}\xspace}  
\newcommand{\Tgammap}[2]{\ensuremath{\Tgamma_{#1#2}}\xspace}  
\newcommand{\Tomegap}[2]{\ensuremath{\Tomega_{#1#2}}\xspace}  

\newcommand{\gtau}{\ensuremath{\mathbf g^{\tau}}\xspace}
\newcommand{\galpha}{\ensuremath{\mathbf g^{\alpha}}\xspace}
\newcommand{\gbeta}{\ensuremath{\mathbf g^{\beta}}\xspace}
\newcommand{\ggamma}{\ensuremath{\mathbf g^{\gamma}}\xspace}
\newcommand{\gomega}{\ensuremath{\mathbf g^{\omega}}\xspace}

\newcommand{\galphap}[1]{\ensuremath{\galpha_{#1}}\xspace}
\newcommand{\gbetap}[1] {\ensuremath{\gbeta_{#1}}\xspace}
\newcommand{\ggammap}[1]{\ensuremath{\ggamma_{#1}}\xspace}
\newcommand{\gomegap}[1]{\ensuremath{\gomega_{#1}}\xspace}

\newcommand{\vtau}{\ensuremath{\mathbf v^{\tau}}\xspace}
\newcommand{\valpha}{\ensuremath{\mathbf v^{\alpha}}\xspace}
\newcommand{\vbeta}{\ensuremath{\mathbf v^{\beta}}\xspace}
\newcommand{\vgamma}{\ensuremath{\mathbf v^{\gamma}}\xspace}
\newcommand{\vomega}{\ensuremath{\mathbf v^{\omega}}\xspace}

\newcommand{\valphap}[1]{\ensuremath{\valpha_{#1}}\xspace}
\newcommand{\vbetap}[1] {\ensuremath{\vbeta_{#1}}\xspace}
\newcommand{\vgammap}[1]{\ensuremath{\vgamma_{#1}}\xspace}
\newcommand{\vomegap}[1]{\ensuremath{\vomega_{#1}}\xspace}

\newcommand{\htau}{\ensuremath{\mathbf h^{\tau}}\xspace}
\newcommand{\halpha}{\ensuremath{\mathbf h^{\alpha}}\xspace}
\newcommand{\hbeta}{\ensuremath{\mathbf h^{\beta}}\xspace}
\newcommand{\hgamma}{\ensuremath{\mathbf h^{\gamma}}\xspace}
\newcommand{\homega}{\ensuremath{\mathbf h^{\omega}}\xspace}

\newcommand{\halphap}[1]{\ensuremath{\halpha_{#1}}\xspace}
\newcommand{\hbetap}[1] {\ensuremath{\hbeta_{#1}}\xspace}
\newcommand{\hgammap}[1]{\ensuremath{\hgamma_{#1}}\xspace}
\newcommand{\homegap}[1]{\ensuremath{\homega_{#1}}\xspace}

\newcommand{\gext}{\ensuremath{\mathbf g_\text{ext}}\xspace}
\newcommand{\vext}{\ensuremath{\mathbf v_\text{ext}}\xspace}
\newcommand{\gint}{\ensuremath{\mathbf g_\text{int}}\xspace}

\newcommand{\hext}{\ensuremath{\mathbf h_\text{ext}}\xspace}

\newcommand{\Deltah}{\ensuremath{\Delta \mathbf h}\xspace}

\newcommand{\zeromat}{\ensuremath{\mathbf 0}\xspace}

\newcommand{\gzero}{\ensuremath{\mathbf g_0}\xspace}
\newcommand{\gone}{\ensuremath{\mathbf g_1}\xspace}
\newcommand{\gtwo}{\ensuremath{\mathbf g_2}\xspace}
\newcommand{\gthree}{\ensuremath{\mathbf g_3}\xspace}

\newcommand{\ftau}{\ensuremath{\mathbf f^{\tau}}\xspace}

\newcommand{\utau}{\ensuremath{\mathbf u^{\tau}}\xspace}
\newcommand{\wtau}{\ensuremath{\mathbf w^{\tau}}\xspace}

\newcommand{\uout}{\ensuremath{\mathbf u^\text{(out)}}\xspace}
\newcommand{\uin}{\ensuremath{\mathbf u^\text{(in)}}\xspace}

\newcommand{\ukout}{\ensuremath{u_k^\text{(out)}}\xspace}  
\newcommand{\ukin}{\ensuremath{u_k^\text{(in)}}\xspace}

\newsiamremark{remark}{Remark}
\newsiamremark{hypothesis}{Hypothesis}
\crefname{hypothesis}{Hypothesis}{Hypotheses}
\newsiamthm{claim}{Claim}

\headers{Quadtree-Adaptive HPS Method for Elliptic PDEs}{D. Chipman, D. Calhoun, and C. Burstedde}

\title{A Fast Direct Solver for Elliptic PDEs on a Hierarchy of Adaptively Refined Quadtrees}

\author{
  Damyn Chipman
    \thanks{Boise State University, Boise, ID (\email{DamynChipman@boisestate.edu}), Corresponding Author.}
  \and
  Donna Calhoun
    \thanks{Boise State University, Boise, ID (\email{DonnaCalhoun@boisestate.edu}).}
  \and
  Carsten Burstedde
    \thanks{Rheinische Friedrich-Wilhelms-Universität Bonn, Germany (\email{burstedde@ins.uni-bonn.de}).}
}

\usepackage{amsopn}

\ifpdf
\hypersetup{
  pdftitle={A Fast and Adaptive Direct Solver for Elliptic PDEs on a Hierarchy of aptively Refined Quadtrees},  pdfauthor={D. Chipman, D. Calhoun, C. Burstedde}}
\fi

\begin{document}

\maketitle


\begin{abstract}
We describe a fast, direct solver for elliptic partial differential equations on a two-dimensional hierarchy of adaptively refined, Cartesian meshes.  Our solver, inspired by the Hierarchical Poincaré-Steklov (HPS) method introduced by Gillman and Martinsson (SIAM J. Sci. Comput., 2014) uses fast solvers on locally uniform Cartesian patches stored in the leaves of a quadtree and is the first such solver that works directly with the adaptive quadtree mesh managed using the grid management library \pforest (C. Burstedde, L. Wilcox, O. Ghattas, SIAM J. Sci. Comput. 2011). Within each Cartesian patch, stored in leaves of the quadtree, we use a second order finite volume discretization on cell-centered meshes.  Key contributions of our algorithm include 4-to-1 merge and split implementations for the HPS build stage and solve stage, respectively.  We demonstrate our solver on Poisson and Helmholtz problems with a mesh adapted to the high local curvature of the right-hand side.
\end{abstract}

\begin{keywords}
Direct Methods, Elliptic Partial Differential Equations, Adaptive Mesh Refinement, Quadtree/Octree Codes
\end{keywords}


\input{sections/introduction}
\input{sections/math-theory}
\input{sections/implementation-details}
\input{sections/results}
\input{sections/conclusion}

\section{Acknowledgements}

The authors wish to thank the support of the National Science Foundation support from grant NSF-DMS \#1819257, and the support of the Boise State University Research Computing facilities.

\bibliographystyle{siamplain}
\bibliography{main}

\end{document}

%% file: sections/introduction.tex
\section{Introduction}
\label{sec:intro}

Elliptic partial differential equations (PDEs) arise in fluid modeling, electromagnetism, astrophysics, heat transfer, and many other scientific and engineering applications. Solving elliptic PDEs numerically can place significant computational demands on scientific codes, and so development of fast, efficient solvers continues to be an active area of research.  Efficient solvers based on finite difference or finite element methods will typically take advantage of the sparsity of the underlying linear systems and any mesh regularity.  Efficient solvers have been developed (UMFPACK \cite{davis2004algorithm}, FISHPACK \cite{swarztrauber1999fishpack}) for uniform Cartesian meshes, since for these problems, the resulting linear system has a regular sparsity pattern that can be exploited.

For many applications, uniform Cartesian meshes are prohibitively expensive and so local mesh adaptivity can be used to resolve detailed solution features.  One common approach to Cartesian mesh adaptivity is the patch-based method, originally developed by Berger, Oligier and later Colella. Software libraries implementing methods related to this patch-based approach include Chombo \cite{colella2009chombo}, \amrex \cite{zhang2019amrex}, and SAMRAI \cite{hornung2006managing}. In a patch-based approach, the computational mesh is defined as a union of overlapping rectangular patches, with finer patches placed on top of coarser patches in regions where resolution for  detailed solution features is needed.  These rectangular patches are of arbitrary size, but align with mesh coordinates of the underlying coarser mesh. Proper nesting rules ensure that every patch is completely contained within coarser patches at the same coarse resolution.

A second approach to Cartesian mesh adaptivity, and the approach used in this paper, is to construct a composite mesh by filling leaves of an adaptive quadtree with non-overlapping locally Cartesian grids of the same size. Quadrants in the quadtree are refined by subdividing the quadrant into four smaller quadrants, each occupying one quarter of the space of the larger quadrant. Similarly, coarsening occurs when four finer quadrants are replaced by a single coarser quadrant.  Every quadrant in the quadtree layout contains a mesh of the same size, but since each quadrant occupies a space determined from their level in the quadtree, the effective resolution of each grid is finer or coarser, depending on their level in the quadtree.  Adaptive, Cartesian software libraries using a tree-based approach include PARAMESH \cite{globisch1995parmesh}, FLASH-X \cite{dubey2022flash}, \pforest \cite{burstedde2011p4est,burstedde2020parallel}, and ForestClaw \cite{calhoun2017forestclaw}.

Solving elliptic problems on an adaptive hierarchy of meshes introduces  technical challenges that are not present with uniform Cartesian meshes. Methods that require matrix assembly are much more difficult to use on adaptive meshes, since row entries corresponding to discretizations at boundaries between coarse and fine meshes are based on non-trivial stencils.  The Hypre library \cite{falgout2002hypre} provides some tools for matrix assembly, but these tools are not immediately useful in adaptive mesh case. The deal.II finite element library uses quadtree/octree mesh adaptivity and performs linear solves with full matrix assembly \cite{bangerth2011algorithms, bangerth2007deal}. A more common approach is to use matrix-free methods such as multigrid or Krylov methods.  Multigrid may be particularly well suited for adaptive meshes, since the adaptive levels are automatically built into the mesh hierarchy.  \amrex, for example, makes extensive use of multigrid for their solvers \cite{zhang2019amrex}.  However, the performance of iterative solvers is largely problem dependent and also affected by irregular stencils at boundaries between coarse and fine meshes.

In \cite{gillman2014direct}, Gillman and Martinsson describe a direct solver for elliptic problems that is particularly well suited to composite quadtree meshes. In their approach, a matrix-free factorization of the linear system is constructed by merging leaf level quadrants in a uniformly refined composite quadtree mesh up to a root node of the tree. Then in a solve stage, the boundary condition data imposed on the physical domain is propagated back down the tree, where Dirichlet problems on each leaf patch can be solved using fast solvers. This approach does not require any matrix assembly on the composite mesh and can be used with any uniform grid solver at the leaf level. In their original work, Gillman and Martinsson use high-order spectral collocation methods and low-rank optimization to achieve $\mathcal O(N)$ factorization complexity \cite{gillman2014direct}.  This method has come to be known as the \HPS (HPS) method \cite{martinsson2015hierarchical}. Advancements to the HPS method include 3D and parallel implementations \cite{hao2016direct,beams2020parallel}. Adaptive mesh variations of the HPS method have been done, and we build on these ideas \cite{babb2018accelerated, geldermans2019adaptive}. Applications of the HPS method can be found in \cite{fortunato2020ultraspherical}.  More details on the HPS method can be found in Chapters 19--27 of \cite{martinsson2019fast} and a tutorial on the HPS method can be found in \cite{martinsson2015hierarchical}. 

The focus of this paper is on an implementation of the HPS method on an adaptive, finite volume mesh using the state-of-the-art meshing library \pforest to generate the adaptive quadtree mesh \cite{burstedde2011p4est,burstedde2020parallel}.  The \pforest library provides highly efficient data structures and iterators necessary for implementing the HPS method.  A particularly useful feature is that \pforest allows for multiple trees (e.g. a ``forest-of-octrees"), allowing for significant flexibility in the geometry of the mesh. By using software libraries like \pforest for the underlying mesh management and PETSc \cite{balay2024petsc} and FISHPACK \cite{swarztrauber1999fishpack} for low-level solvers, we can build a modular codebase that can be applied to a variety of applications. This implementation is built into the EllipticForest \cite{chipman2023elliptic} software library for solving elliptic PDEs on adaptive meshes.

\subsection{Problem Statement}

We focus on the constant coefficient elliptic problem
\begin{equation}
    \nabla^2 u(x,y) + \lambda u(x,y) = f(x,y), \quad \lambda \ge 0
    \label{eq:elliptic_pde}
\end{equation}
with $(x,y) \in \Omega = [a, b] \times [c,d]$, subject to Dirichlet boundary conditions
\begin{equation}
    u(x,y) = g(x,y),\ \ (x,y) \in \Gamma = \partial \Omega.
\end{equation}
For this paper, we assume that $b-a=d-c$ so that we can describe our mesh using a single quadtree and square mesh cells.  In general though, this is not a limitation of meshes generated using the multi-block ``forest-of-octrees'' capabilities of \pforest. In addition, we may use mappings between square computational domains to moderately rectangular cells for more complex meshes. The primary contributions of this paper are an HPS solver for the constant coefficient elliptic problem on an adaptively refined \pforest mesh.

In \refsec{sec:quadtree}, we describe the quadtree-adaptive HPS method, including the leaf-level computations, merging and splitting algorithms, and the handling of coarse-fine interfaces. This section also includes a comparison between our 4-to-1 approach contrasted to the 2-to-1 approach presented in \cite{gillman2014direct}. \refsec{sec:adaptivity} describes implementation details such as the data structures and functions that wrap the \pforest mesh backend. Finally, in \refsec{sec:results} we provide results of numerical experiments solving different types of elliptic PDEs, including a convergence analysis and timing results for the use of the quadtree-adaptive HPS method on adaptive meshes.

%% file: sections/math-theory.tex
\section{Components of the HPS algorithm: Merging, Splitting and Leaf Computations}
\label{sec:quadtree}
The HPS method applied to \refeq{eq:elliptic_pde} includes a factorization stage and a solve stage. For problems in which for $f(x,y)$ is non-zero, an additional ``upwards'' stage is required for a particular right hand side.  In our algorithm, the build and upwards stages are a successive application of a 4-to-1 merge algorithm while the solve stage is a an application of a 1-to-4 split algorithm.  In this section, we describe these merging and splitting operations.

\ignore{We begin by detailing the computation required on any leaf level patches. We then outline the algorithms for the merging and splitting of families of patches. Next, we demonstrate the HPS method with the build, upwards, and solve stages. }

\subsection{Leaf Level Computations}
\label{sub:leaf_level_computations}

\begin{figure}
    \centering
    \begin{tabular}{ccc}
        \begin{subfigure}[t]{0.3\textwidth}
            \centering
            \includegraphics[width=\textwidth, clip=true, trim={100 150 100 150}]{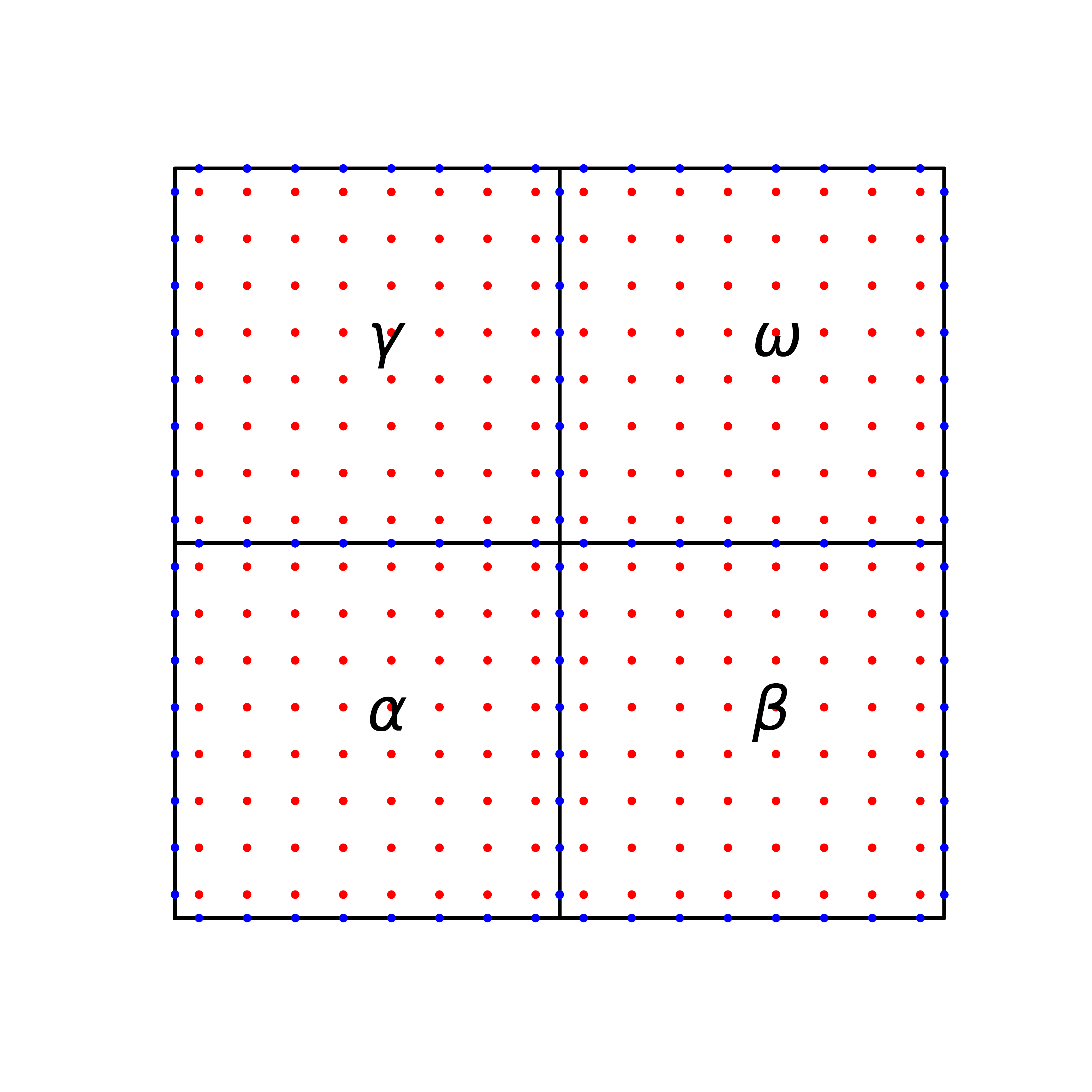}
            \label{subfig:4_patches_with_grid}
        \end{subfigure}
        &
        \begin{subfigure}[t]{0.3\textwidth}
            \centering
            \includegraphics[width=\textwidth, clip=true, trim={100 150 100 150}]{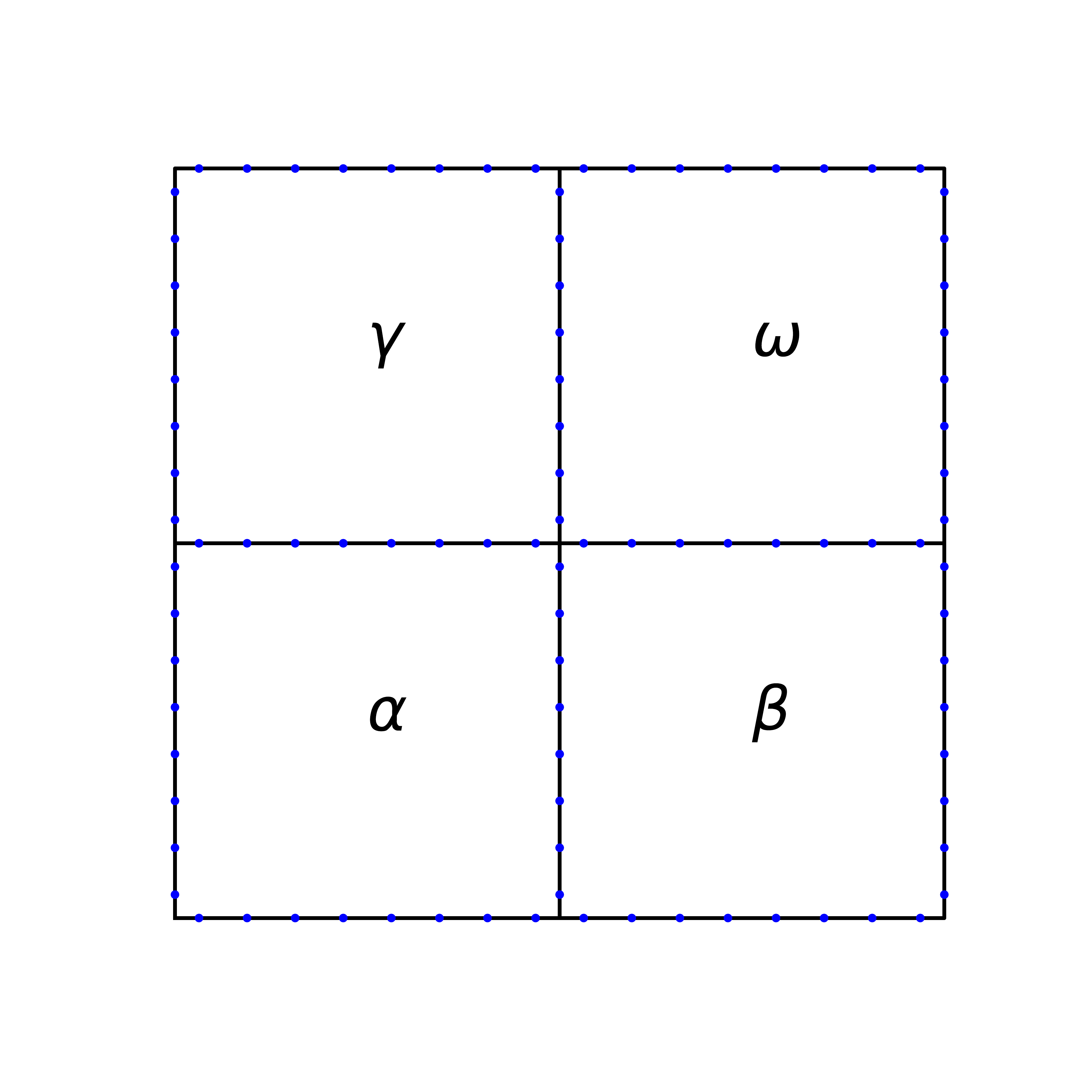}
            \label{subfig:4_patches}
        \end{subfigure}
        &
        \begin{subfigure}[t]{0.3\textwidth}
            \centering
            \includegraphics[width=\textwidth, clip=true, trim={100 150 100 150}]{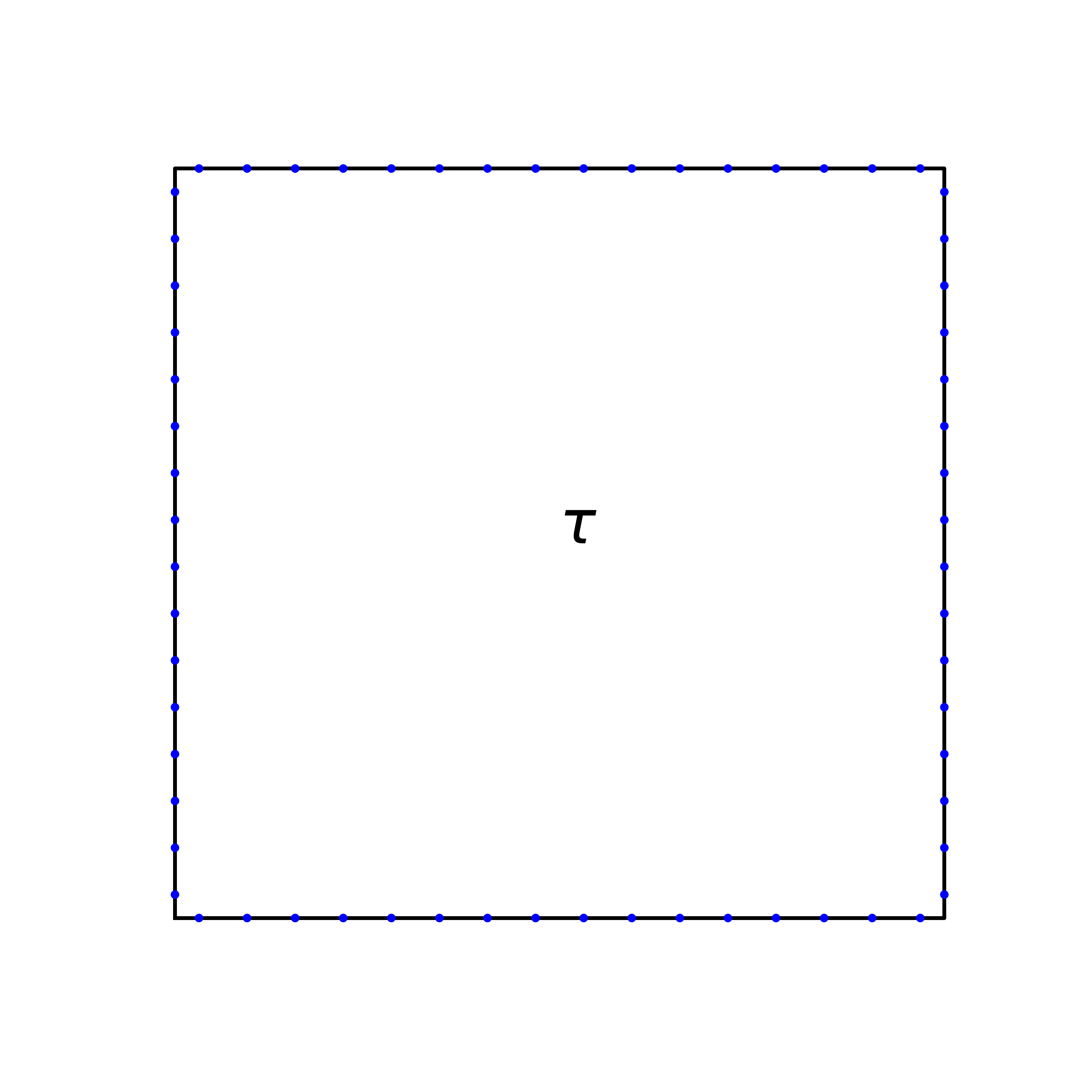}
            \label{subfig:parent_patch}
        \end{subfigure}
    \end{tabular}
    \caption{The 4-to-1 merge process: (left) the four children patches with their local grid, (middle) the four children share internal and external boundaries with $\tau$, (right) the merged parent patch with data on the exterior of $\tau$.}
    \label{fig:4_to_1_patches}
\end{figure}

In what follows, we assume that each leaf-level quadrant in the quadtree mesh stores a uniform Cartesian grid with $M \times M$ mesh cells.  We refer to this local Cartesian mesh, along with its quadrant as a {\em patch}.  In \reffig{fig:4_to_1_patches} (left), we show four patches in a larger composite domain. Cell-centered finite difference schemes are particularly convenient for adaptively refined Cartesian meshes, since the cell-centered values are not duplicated on adjacent patches.

At the leaf-level, the build stage of the HPS methods requires the solution to a Poisson problem \eqn{elliptic_pde} discretized on an $M \times M$ cell-centered grid as
\begin{equation}
\Atau \utau = \Bctau \gtau + \ftau
\label{eq:patch_poisson_problem}
\end{equation}
where $\Bctau$ spreads Dirichlet boundary data $\gtau \in \real^{4M}$ to the grid and \ftau is the right hand side function $f(x,y)$ from \eqn{elliptic_pde} evaluated at cell centers.  We write the solution to this problem in terms of a {\em homogeneous operator} \Lhtau and an inhomogeneous operator \Linhtau which represent efficient solvers on the uniform Cartesian patch. The solution is then given by 
\begin{equation}
\utau = \Lhtau \gtau + \Linhtau \ftau
\label{eq:patch_leaf_solve}
\end{equation}
Where the context is clear, we also view operators \Lhtau and \Linhtau, boundary data \gtau and right hand side data \ftau as matrices and vectors of the appropriate sizes. The components of \utau are given by $u^\tau_{ij}$, for $i,j = 1,2,\hdots M$.  

Using these leaf-level solvers, we can now build a leaf-level DtN operator \Ttau.  In the description below, we illustrate the HPS method on the complete Poisson problem, assuming a non-zero right hand side.  We remark below, though, that if one anticipates solving for multiple right-hand sides, the build stage can be separated into a factorization stage and an "upwards" stage to handle inhomogeneous data.  

\ignore{\donna{Not sure where to put this.}  The general HPS algorithm does not specify a particular local mesh discretization, and several have been proposed in the literature.  In \cite{gillman2014direct}, the authors use a spectral collocation method and in \cite{fortunato2020ultraspherical}, a Chebyshev tensor product grid for high order finite elements was used. For compatibility with the finite volume code ForestClaw \cite{calhoun2017forestclaw}, we use a second order finite volume discretization of \refeq{eq:elliptic_pde}.}

\subsubsection{\DtN (DtN) operator}
Given a cell-centered grid solution $\utau = \Lhtau \gtau + \Linhtau \ftau$, we define ``interior" grid values $\uin \in \real^{4M}$ as those grid solution values \utau for which $i = 1$, $i = M$, $j = 1$ or $j = M$.  By analogy, components of \uout are grid solution values for which $i=0$ or $i = M+1$ or $j = 0$ or $j = M+1$.     

On a cell-centered finite volume grid, the known Dirichlet boundary values are not collocated with components of the grid solution so we discretize the Dirichlet boundary condition at the midpoint between \ukin and \ukout as
\begin{equation}
g_k^\tau = \frac{\ukout + \ukin}{2}, \quad k = 1,\hdots , 4M
\end{equation}
where $g_k$, $\ukin$ and $\ukout$ are the components of \gtau, \uin and \uout respectively.  

By analogy, we discretize Neumann data on the boundary of the patch as
\begin{equation}
v_k = \frac{\ukout - \ukin}{h}
\end{equation}
where $h$ is the mesh cell width (assumed to be uniform in both x- and y- directions).

Eliminating $\ukout$  between the two expressions, we  can express the Neumann data $v_k$ in terms of Dirichlet data $g_k$ and the interior value $\ukin$ as
\begin{equation}
v_k = \frac{2}{h}(g_k - \ukin).
\end{equation}
The vector of Neumann data on patch $\tau$ can then be expressed as 
\begin{equation}
\vtau = \frac{2}{h}(\gtau - \uin).
\label{eq:vtau}
\end{equation}   

Introducing an operator (or a matrix of the appropriate size) \Gop that selects \uin from the grid solution data \utau, we write $\uin = \Gop \utau$ so that
\begin{equation}
\vtau = \frac{2}{h}(\gtau - \Gop (\Lhtau \gtau + \Linhtau \ftau)) = 
\frac{2}{h}(I - \Gop \Lhtau)\gtau - \frac{h}{2}\Gop \Linhtau \ftau
\end{equation}
From this, we define the leaf level \DtN mapping $\Ttau \in \real^{4M \times 4M}$ as
\begin{equation}
\Ttau \equiv \frac{2}{h}(I - \Gop\Lhtau).
\label{eq:map_D2N}
\end{equation}
with inhomogenenous component
\begin{equation}
\htau \equiv -\frac{2}{h}\Gop\Linhtau \ftau. 
\end{equation}
Neumann data on the patch is then computed as
\begin{equation}
\vtau = \Ttau \gtau + \htau.
\label{eq:vtau_with_h}
\end{equation}
We refer \Ttau as a discrete \DtN operator, although the Neumann data \vtau is not the Neumann data that is commonly understood to be the result of applying a DtN map, since \vtau depends not only on the Dirichlet data \gtau, but also on the inhomogeneous data \ftau.  

\subsection{The 4-to-1 Merge Algorithm}
\label{sub:4-to-1merge}
If $\tau$ is not a leaf-level quadrant, then we construct the DtN map \Ttau {\em recursively} by merging mappings from the four children patches $\alpha$, $\beta$, $\gamma$ and $\omega$  partitioning $\tau$.  For the non-leaf, two additional operators $\Stau$ and $\Xtau$ (for non-zero $f(x,y)$) are also constructed.  Prior to this merge process, it is assumed that each patch has computed a DtN mapping \Ti, for $i=\alpha$, $\beta$, $\gamma$, $\omega$. 

\subsubsection{Merging DtN mappings} We partition sibling operators \Ti, $i=\alpha$, $\beta$, $\gamma$, $\omega$ according to how Dirichlet data on an edge of patch $i$ are mapped to an edge on $i$ containing Neumann data.  For example, the sub-matrix $\Talphap{\tau}{\gamma}$ maps Dirichlet data on the boundary shared between sibling patch $\alpha$  and parent patch $\tau$ to Neumann data on the boundary shared between patch $\alpha$ and sibling patch $\gamma$. In a similar manner, the Dirichlet data on each of the sibling patches $i$ into $\galphap{\tau}$, $\galphap{\gamma}$ and $\galphap{\beta}$, where the subscripts indicate data on edges that patch $\alpha$ shares with patches $\tau$, $\gamma$ and $\beta$. Neumann data $\valpha$ and \halpha are partitioned in an analogous manner. 

With these partitions, we can write a set of equations for each of the four sibling DtN mappings as
\begin{equation}
\begin{aligned}
\valpha = \Talpha \galpha + \halpha \quad & \Longrightarrow  & \quad 
    \begin{bmatrix}
    \valphap{\tau} \\
    \valphap{\gamma} \\
    \valphap{\beta} \\
    \end{bmatrix}
    =
    \begin{bmatrix}
    \Talphap{\tau}{\tau}   &\Talphap{\gamma}{\tau}   & \Talphap{\beta}{\tau} \\
    \Talphap{\tau}{\gamma} &\Talphap{\gamma}{\gamma} & \Talphap{\beta}{\gamma} \\
    \Talphap{\tau}{\beta}  &\Talphap{\gamma}{\beta}  & \Talphap{\beta}{\beta} \\
    \end{bmatrix}
    \begin{bmatrix}
    \galphap{\tau} \\
    \galphap{\gamma} \\
    \galphap{\beta} \\
    \end{bmatrix} + 
    \begin{bmatrix}
    \halphap{\tau} \\
    \halphap{\gamma} \\
    \halphap{\beta} \\
    \end{bmatrix}\\
&& \\
\vbeta = \Tbeta \gbeta + \hbeta \quad & \Longrightarrow  & \quad
    \begin{bmatrix}
        \vbetap{\tau} \\
        \vbetap{\omega} \\
        \vbetap{\alpha} \\
    \end{bmatrix}
    =
    \begin{bmatrix}
        \Tbetap{\tau}{\tau}   & \Tbetap{\omega}{\tau}   & \Tbetap{\alpha}{\tau} \\
        \Tbetap{\tau}{\omega} & \Tbetap{\omega}{\omega} & \Tbetap{\alpha}{\omega} \\
        \Tbetap{\tau}{\alpha} & \Tbetap{\omega}{\alpha} & \Tbetap{\alpha}{\alpha} \\
    \end{bmatrix}
    \begin{bmatrix}
        \gbetap{\tau} \\
        \gbetap{\omega} \\
        \gbetap{\alpha} \\
    \end{bmatrix} + 
    \begin{bmatrix}
        \hbetap{\tau} \\
        \hbetap{\omega} \\
        \hbetap{\alpha}
    \end{bmatrix}\\
&& \\
\vgamma = \Tgamma \ggamma + \hgamma \quad & \Longrightarrow  & \quad
    \begin{bmatrix}
        \vgammap{\tau} \\
        \vgammap{\alpha} \\
        \vgammap{\omega} \\
    \end{bmatrix}
    =
    \begin{bmatrix}
        \Tgammap{\tau}{\tau}    & \Tgammap{\alpha}{\tau}   & \Tgammap{\omega}{\tau} \\
        \Tgammap{\tau}{\alpha}  & \Tgammap{\alpha}{\alpha} & \Tgammap{\omega}{\alpha} \\
        \Tgammap{\tau}{\omega}  & \Tgammap{\alpha}{\omega} & \Tgammap{\omega}{\omega} \\
    \end{bmatrix}
    \begin{bmatrix}
        \ggammap{\tau} \\
        \ggammap{\alpha} \\
        \ggammap{\omega} \\
    \end{bmatrix} + 
    \begin{bmatrix}
        \hgammap{\tau} \\
        \hgammap{\alpha} \\
        \hgammap{\omega} \\
    \end{bmatrix}\\    
&& \\    
\vomega = \Tomega \gomega + \homega \quad & \Longrightarrow  & \quad
    \begin{bmatrix}
        \vomegap{\tau} \\
        \vomegap{\beta} \\
        \vomegap{\gamma} \\
    \end{bmatrix}
    =
    \begin{bmatrix}
        \Tomegap{\tau}{\tau}   & \Tomegap{\beta}{\tau}   & \Tomegap{\gamma}{\tau} \\
        \Tomegap{\tau}{\beta}  & \Tomegap{\beta}{\beta}  & \Tomegap{\gamma}{\beta} \\
        \Tomegap{\tau}{\gamma} & \Tomegap{\beta}{\gamma} & \Tomegap{\gamma}{\gamma} \\
    \end{bmatrix}
    \begin{bmatrix}
        \gomegap{\tau} \\
        \gomegap{\beta} \\
        \gomegap{\gamma} \\
    \end{bmatrix} +
    \begin{bmatrix}
        \gomegap{\tau} \\
        \gomegap{\beta} \\
        \gomegap{\gamma} \\
    \end{bmatrix} \\
\end{aligned}.
\label{eq:linalg_one}
\end{equation}
Taking the first equation from each of the four sets of equations in \eqn{linalg_one}, we write a system of equations for the Dirichlet data to get
\begin{equation}
A\gext + B \gint + \hext = \vext
\label{eq:sys_one}
\end{equation}
where the Dirichlet data has been partitioned into {\em exterior} data and {\em interior} data.  The exterior data is data on the boundary of the parent patch $\tau$ is given by
\begin{equation}
\gext \equiv \left[\galphap{\tau}, \gbetap{\tau} , \ggammap{\tau}, \gomegap{\tau}\right]^T
\end{equation}
To define the interior data, we use the fact that our solution is continuous across boundaries shared by sibling grids and enumerate data on the four shared boundaries as
\begin{equation}
\begin{aligned}
\gzero \equiv \galphap{\gamma} & = \ggammap{\alpha}\\
\gone \equiv \gbetap{\omega} & = \gomegap{\beta}  \\
\gtwo \equiv \galphap{\beta} & = \gbetap{\alpha} \\
\gthree \equiv \ggammap{\omega} & = \gomegap{\gamma}.
\end{aligned}
\end{equation}
With the Dirichlet data on the shared boundary uniquely defined, the interior partition of the Dirichlet data is given as
\begin{equation}
\gint \equiv \left[\gzero, \gone, \gtwo, \gthree\right]^T
\end{equation}
The Neumann data \vtau and \htau is partitioned analogously. 

The block matrices $A$ and $B$ are defined as
\begin{equation}
A = \begin{bmatrix}
\Talphap{\tau}{\tau} & \zeromat & \zeromat & \zeromat \\
\zeromat & \Tbetap{\tau}{\tau} & \zeromat & \zeromat  \\
\zeromat & \zeromat & \Tgammap{\tau}{\tau} & \zeromat  \\
\zeromat& \zeromat & \zeromat & \Tomegap{\tau}{\tau} \\
\end{bmatrix}, \qquad
B = \begin{bmatrix}
\Talphap{\gamma}{\tau} & \zeromat              & \Talphap{\beta}{\tau} & \zeromat \\
\zeromat               & \Tbetap{\omega}{\tau} & \Tbetap{\alpha}{\tau} & \zeromat \\
\Tgammap{\alpha}{\tau} & \zeromat              & \zeromat              & \Tgammap{\omega}{\tau} \\
\zeromat               & \Tomegap{\beta}{\tau} & \zeromat              & \Tomegap{\gamma}{\tau}
\end{bmatrix}
\label{eq:matrix_AB}
\end{equation}

To obtain a second set of equations, we impose a continuity condition on the normal derivative across edges shared between sibling patches and get
\begin{equation}
\begin{aligned}
\valphap{\gamma} + \vgammap{\alpha} & = \zeromat \\
\vbetap{\omega} + \vomegap{\beta} & = \zeromat \\
\valphap{\beta}  + \vbetap{\alpha} & = \zeromat\\
\vgammap{\omega} + \vomegap{\gamma} & = \zeromat.
\end{aligned}
\end{equation}
We organize the resulting four equations as
\begin{equation}
C \gext + D\gint  + \Deltah = \zeromat
\label{eq:sys_two}
\end{equation}
where
\begin{equation}
C = 
\begin{bmatrix}
\Talphap{\tau}{\gamma} & \zeromat               & \Tgammap{\tau}{\alpha} & \zeromat \\
\zeromat               & \Tbetap{\tau}{\omega}  & \zeromat                & \Tomegap{\tau}{\beta} \\ 
\Talphap{\tau}{\beta}  & \Tbetap{\tau}{\alpha} & \zeromat                & \zeromat \\
\zeromat               & \zeromat               & \Tgammap{\tau}{\omega}  & \Tomegap{\tau}{\gamma},
\end{bmatrix}
\end{equation}
\begin{equation}
D = \begin{bmatrix}
\Talphap{\gamma}{\gamma} + \Tgammap{\alpha}{\alpha} 
& \zeromat 
& \Talphap{\beta}{\gamma} 
& \Tgammap{\omega}{\alpha} \\
\zeromat 
& \Tbetap{\omega}{\omega} + \Tomegap{\beta}{\beta} 
& \Tbetap{\alpha}{\omega} 
& \Tomegap{\gamma}{\beta} \\
\Talphap{\gamma}{\beta} 
& \Tbetap{\omega}{\alpha} 
& \Talphap{\beta}{\beta}+ \Tbetap{\alpha}{\alpha} 
& \zeromat \\
\Tgammap{\alpha}{\omega} 
& \Tomegap{\beta}{\gamma} 
& \zeromat 
& \Tgammap{\omega}{\omega} + \Tomegap{\gamma}{\gamma}
\end{bmatrix}
\label{eq:matrix_CD}.
\end{equation}
and 
\begin{equation}
\Deltah = 
\begin{bmatrix}
\halphap{\gamma} + \hgammap{\alpha} \\
\hbetap{\omega} + \homegap{\beta} \\
\halphap{\beta} + \hbetap{\alpha} \\
\hgammap{\omega} + \homegap{\gamma}
\end{bmatrix}.
\end{equation}

Combining \eqn{sys_one} and \eqn{sys_two}, we have
\begin{equation}
\begin{aligned}
A\gext + B\gint + \hext & = \vext \\
C\gext + D\gint + \Deltah & = \zeromat
\end{aligned}
\end{equation}
Writing this system as an augmented system and applying a single step of block Gaussian elimination, we get
\begin{equation}
\left[
\begin{array}{cc|c}
A & B & \vext - \hext\\
C & D & -\Deltah\\
\end{array}\right]
\;
\Longrightarrow 
\;
\left[
\begin{array}{cc|c}
A - BD^{-1}C & 0 & \vext - \hext + B D^{-1}\Deltah\\
D^{-1}C      & I & -D^{-1}\Deltah
\end{array}\right]
\label{eq:reduced_system}
\end{equation}
We express the first row of the reduced system in equation form as 
\begin{equation}
\left(A - BD^{-1}C\right) \gext = \vext - \hext + B D^{-1}\Deltah
\label{eq:complete_ops}
\end{equation}
From this, we can identify the {\em merged} DtN operator 
\begin{equation}
\Ttau \equiv \mathbf{A} - \mathbf{B} \mathbf{D}^{-1} \mathbf{C} 
\label{eq:dtn_merge}
\end{equation}
as the operator that maps exterior Dirichlet data to exterior Neumann data, with merged inhomogeneous Neumann data given as
\begin{equation}
\htau \equiv -\hext + B D^{-1}\Deltah. 
\end{equation}
With these choices of \Ttau and \htau  we recover the expression from \eqn{vtau}.

The second row of the reduced system \eqn{reduced_system}, given by 
\begin{equation}
\gint = -D^{-1}C \gext -D^{-1} \Deltah,
\end{equation}
shows us how to recover interior Dirichlet data \gint  from exterior data \gext.  
From this, we introduce the {\em solve} operator \Stau, given by 
\begin{equation}
\Stau \equiv -D^{-1}C
\end{equation}

Introducing  the operator
\begin{equation}
\Xtau \equiv -D^{-1}
\label{eq:xtau}
\end{equation}
we can write \Ttau and  \Stau in convenient form as
\begin{equation}
\begin{aligned}
\Stau & = \Xtau C \\
\Ttau & = A + B\Stau
\label{eq:stau_ttau}
\end{aligned}
\end{equation}
To compute the merged inhomogeneous term, we first compute
\begin{equation}
\wtau \equiv \Xtau\Deltah
\label{eq:wtau}
\end{equation}
so that
\begin{equation}
\htau = -\hext - B\wtau
\label{eq:htau}
\end{equation}

Once the build stage is complete, the solve state maps Dirichlet data on parent quadrants to Dirichlet data on child quadrants via the solve operator \Stau, as
\begin{equation}
\gint = \Stau \gext + \wtau.
\label{eq:solve_eqn}
\end{equation}

The formalism presented in terms of merged components \Xtau, \Stau, \Ttau and \htau is identical to that presented in \cite{martinsson2015hierarchical} for the horizontal merge of two leaf boxes.

The build stage of HPS method is described in \Alg{build_merge_uniform}
\begin{algorithm}[H]
    \caption{Build stage on a uniformly refined quadtree mesh}
    \begin{algorithmic}[0]
        \For{$\tau = 0,1,\dots$} \Comment{Iterate over quadrants in uniform quadtree}
            \If{$\tau$ is a leaf-level patch}
                \State Build and store leaf-level DtN map \Ttau and inhomogeneous data \htau
            \Else \Comment{$\tau$ has four child patches $i$, $i=\alpha$, $\beta$, $\gamma$, $\omega$.}
                \State Solve $D\Stau = -C$ to get operator \Stau. 
                \State Build operator \Ttau using \eqn{stau_ttau}.      
                \State Solve $D\wtau = \Deltah$ using to get \wtau.  
                \State Build merged inhomogeneous vector \htau using \eqn{htau}.
                \State Store \Stau, \Ttau, \htau and \wtau with quadrant $\tau$.
                \State Delete operators \Talpha, \Tbeta, \Tgamma, \Tomega.
                \State Delete inhomogeneous data \halpha, \hbeta, \hgamma and \homega.
            \EndIf
        \EndFor
    \end{algorithmic}
    \label{alg:build_merge_uniform}
\end{algorithm}
At the end of this build stage, we will have a single DtN operator \Ttau and single vector \htau for the root patch.  At all levels, though, we will have a solve operator \Stau and inhomogeneous vector \wtau.

A solve stage of the HPS method starts with data $g(x,y)$ defined on the boundary of the domain and successively ``splits" the domain by mapping this exterior Dirichlet data to interior data.  At the leaf level, we solve the Poisson problem \eqn{patch_poisson_problem}.  The solve stage algorithm is described in \Alg{solve_uniform}.
\begin{algorithm}[H]
    \caption{Solve stage on a uniformly refined quadtree mesh}
    \begin{algorithmic}[0]
        \For{$\tau = N,N-1,\hdots,0$} \Comment{Iterate over quadrants in uniform quadtree}
            \If{$\tau$ is a leaf-level patch}
                \State Solve patch Poisson problem \eqn{patch_poisson_problem}
            \Else{} \Comment{$\tau$ has four child patches $\alpha$, $\beta$, $\gamma$, $\omega$.}
                \State Use \eqn{solve_eqn} to map Dirichlet \gext on $\tau$ to Dirichlet data \gint on the interior shared child patch boundaries.
            \EndIf
        \EndFor
    \end{algorithmic}
    \label{alg:solve_uniform}
\end{algorithm}

\remark{If the HPS solver is to be used with multiple right hand-sides for the same composite quadtree mesh, the build stage as described above should be separated into a ``factorization" stage that only stores the operator \Stau and a second ``upwards" stage that builds the inhomogeneous data \htau from \wtau for each patch. Additionally, since \Xtau is needed to compute \wtau, and $B$ is needed to compute \htau from \wtau, one must also store \Xtau and $B$.  Alternatively, one can store \Ttau at each level and rebuild \Xtau, $B$, \wtau and \htau as needed. If the problem in \eqn{elliptic_pde} has zero right-hand side data, any operations involving \htau can be eliminated, since the inhomogeneous Neumann data will be identically zero.}

\remark{For a constant coefficient problem on a uniformly refined mesh, only one operator \Ttau per level needs to be computed and stored. In merge stage, this single operator can be used in place of each child patch. This significantly reduces the time needed in the build stage.}
\label{remark:caching}

\ignore{
\remark{Neumann boundary data $v_k$ supplied along any segment of the domain boundary in \eqn{elliptic_pde} can be easily converted to a Dirichlet data $g_k$ using the formula
\begin{equation}
g_k = \ukin + \frac{h}{2}v_k
\end{equation}
This Dirichlet data is then propagated down to the leaf-level patches using \Alg{solve_uniform}.}}

\remark{Neumann boundary data on edges of the physical domain can be handled using \eqn{vtau_with_h} and  inverting the DtN mapping \Ttau to map Neumann data to Dirichlet data along any edges with Neumann data.  For this, we have only considered the non-singular case and assume that the problem is not purely a Neumann problem.}

\subsection{Coarse-Fine Interfaces}
\label{sub:mesh_adaptivity}

\begin{figure}
    \centering
    \begin{tabular}{ccc}
        \begin{subfigure}[t]{0.3\textwidth}
            \centering
            \includegraphics[width=\textwidth, clip=true, trim={100 150 100 150}]{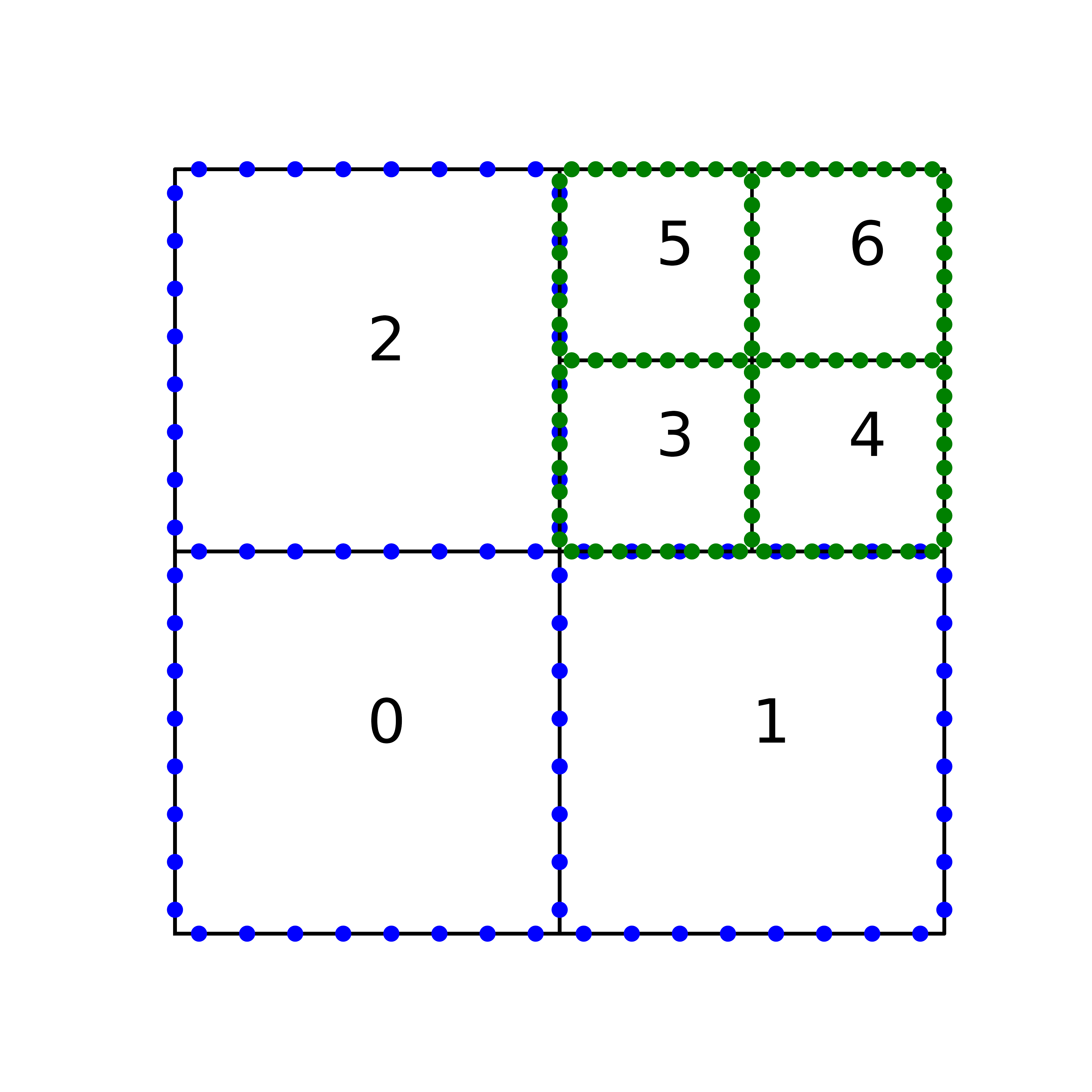}
        \end{subfigure}
        &
        \begin{subfigure}[t]{0.3\textwidth}
            \centering
            \includegraphics[width=\textwidth, clip=true, trim={100 150 100 150}]{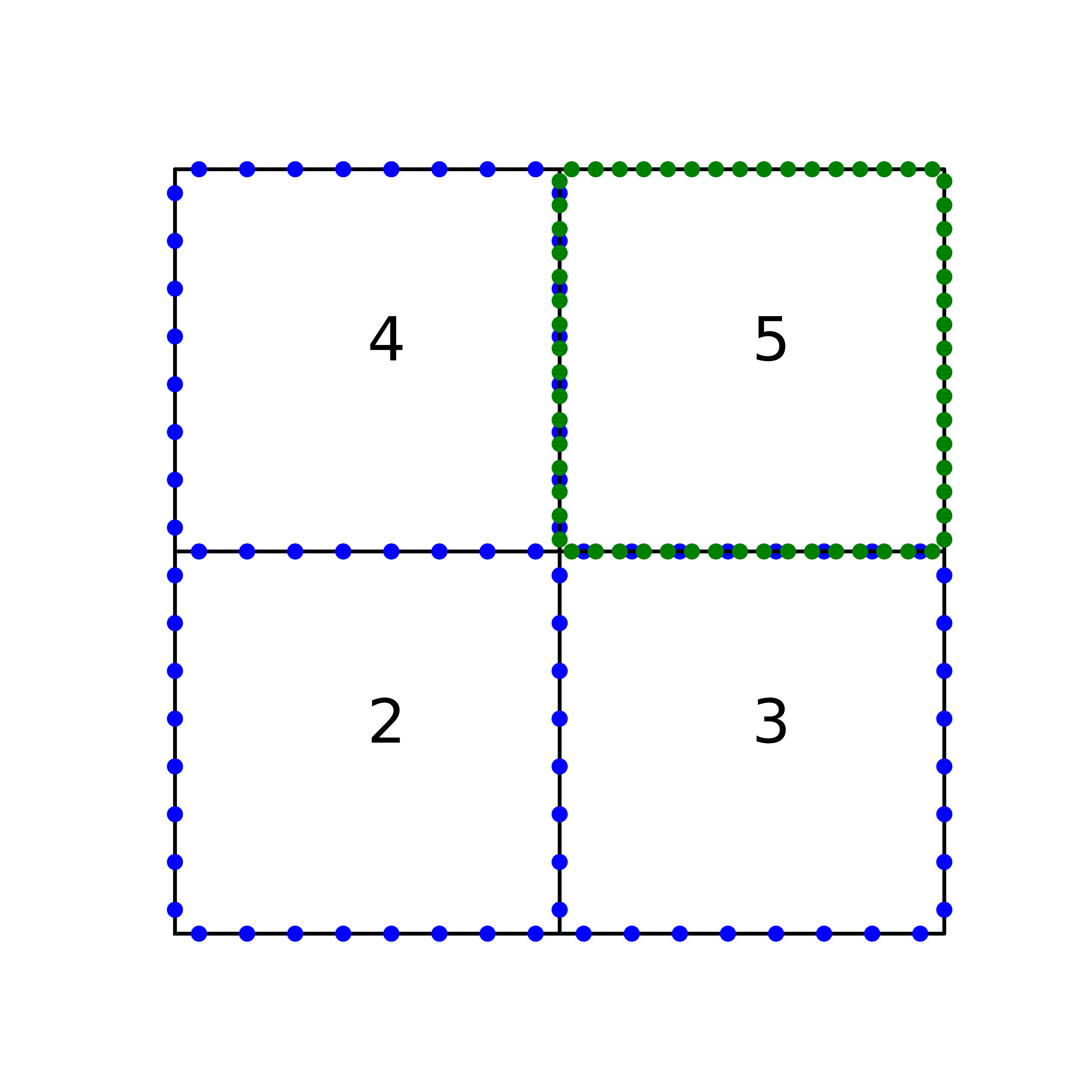}
        \end{subfigure}
        &
        \begin{subfigure}[t]{0.3\textwidth}
            \centering
            \includegraphics[width=\textwidth, clip=true, trim={100 150 100 150}]{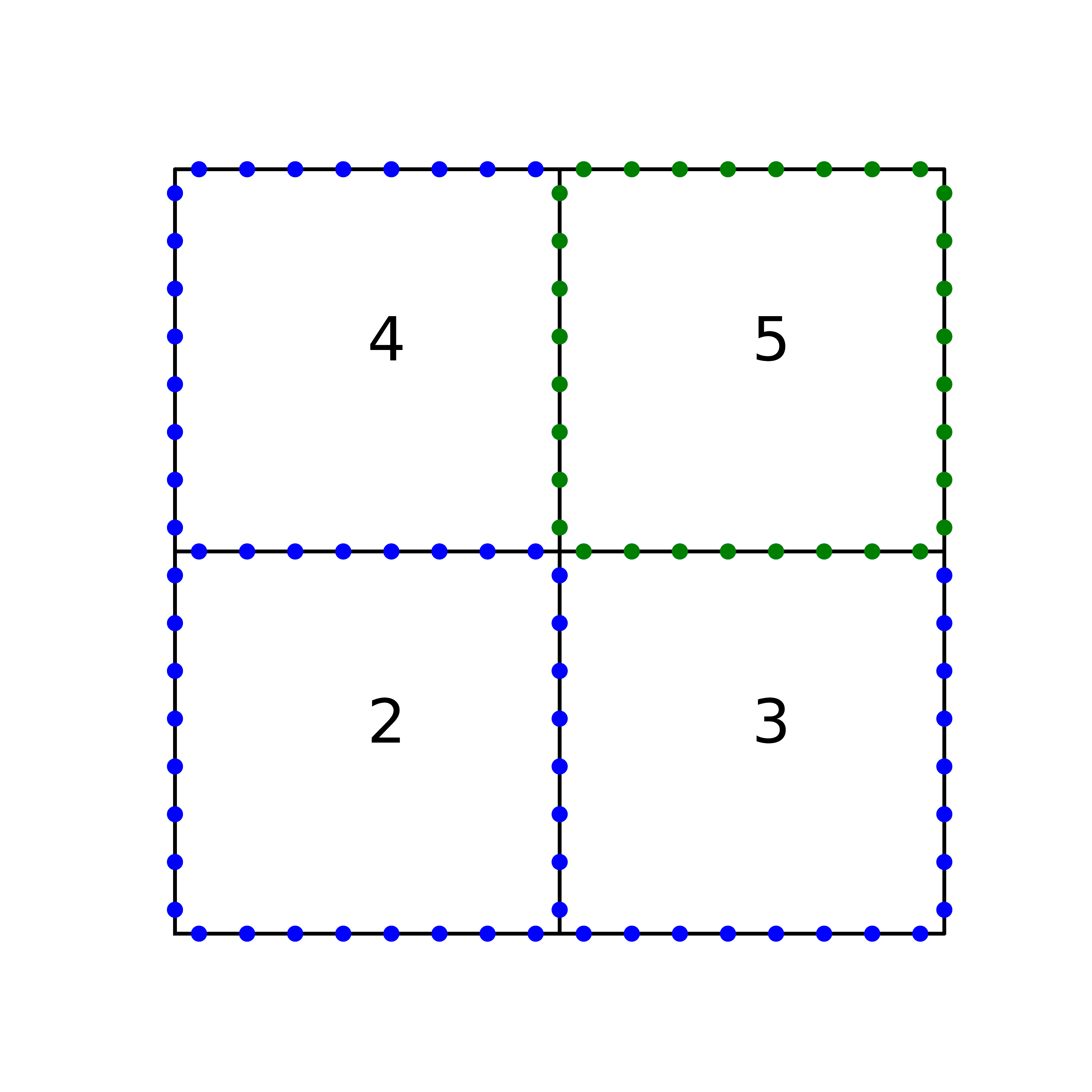}
        \end{subfigure}
    \end{tabular}
    \caption{Example of a coarse-fine interface merge: (left) Patches $4, 5, 6, \&\ 7$ will be merged into parent patch $3$, resulting in a coarse-fine interface. (middle) The data on patch 3 is averaged (coarsened). (right) Merging $0, 1, 2, \&\ 3$ can continue as detailed in \refsec{sub:4-to-1merge}.}
    \label{fig:adaptive_merge}
\end{figure}

The merging and splitting processes we have described so far applies only to four sibling patches occupying a coarser level quadrant. For a fully adaptive quadtree, however, not all quadrants are refined to the same level, resulting in coarse-fine interfaces (see \reffig{fig:adaptive_merge}). To accommodate for the coarse-fine interfaces, we coarsen the data associated with a patch to match the size of it's siblings. The data that needs to be coarsened during the build stage is the DtN matrix \Ttau and the inhomogenenous Neumann data \htau. This approach is the one used when dealing with coarse-fine interfaces in \cite{babb2018accelerated}. Another approach is the one found in \cite{geldermans2019adaptive} where they coarsen the operators along each individual side with a coarse-fine interface, as opposed to the entire patch.

Prior to merging, we check each of the four children patches to be merged for a patch with data more fine than the other children. We do this by checking the size of the associated grid on that patch. If a coarse-fine interface exists, the data on the fine grid patch is coarsened to match its siblings. We form interpolation matrices $\mathbf L_{2,1}$ and $\mathbf L_{1,2}$ that map two sides to one or one side to two, respectively. For our cell-centered finite volume mesh, these matrices either average two data points to one, or interpolate one data point to two. Block diagonal versions of these operators are formed as
\begin{align}
    \mathbf L_{2,1}^{B} &= \texttt{BlockDiag}(\{\mathbf L_{2,1}, \mathbf L_{2,1}, \mathbf L_{2,1}, \mathbf L_{2,1}\}) \\
    \textbf{L}_{1,2}^{B} &= \texttt{BlockDiag}(\{\textbf{L}_{1,2}, \textbf{L}_{1,2}, \textbf{L}_{1,2}, \textbf{L}_{1,2}\})
\end{align}
and coarsening $\textbf{T}$ and $\textbf{h}$ is done by matrix multiplication
\begin{align}
    \textbf{T}^{\tau'} &= \textbf{L}_{2,1}^{B} \textbf{T}^{\tau} \textbf{L}_{1,2}^{B} \\
    \textbf{h}^{\tau'} &= \textbf{L}_{2,1}^{B} \textbf{h}^{\tau}
\end{align}
where the prime indicates a coarsened version of that data.

In practice, it is possible to simply replace the fine version of the DtN matrix $\textbf{T}^{\tau'}$ with a coarsened one that has been pre-computed and cached, as mentioned in Remark \ref{remark:caching}.

A 1-to-4 split with a parent that has data that was coarsened during the 4-to-1 merge (for example, patch 5 from \reffig{fig:adaptive_merge}) needs to be refined. The  data that is refined during the solve stage is the boundary data, \gext in Equation \refeq{eq:solve_eqn}. This is done with the block diagonal interpolation matrices from above:
\begin{align}
    \textbf{g}^{\tau}_{ext} = \textbf{L}_{1,2}^{B} \textbf{g}^{\tau'}_{ext}.
\end{align}

\subsection{Comparison Between 4-to-1 Merging and 2-to-1 Merging}
\label{sub:comparison_between_4t1_and2t1_merging}

Here, we compare the performance and storage requirements for the 4-to-1 merge outlined in this paper against the 2-to-1 merge presented in \cite{gillman2014direct}. In \cite{gillman2014direct}, merging is done between two neighboring patches. To merge a family of four patches, one must do two vertical merges (merge two neighboring patches that lie next to each other in the y-direction) and then one horizontal merge (merge two ``tall-skinny'' patches that lie next to each other in the x-direction). Thus, computing and storing $\textbf{T}^{\tau}$ and $\mathbf{S}^{\tau}$ requires three, 2-to-1 merges: two vertical merges and one horizontal merge.

For both approaches, we assume that a patch has $M$ points per side, resulting in $M^2$ points per patch. We will compare the floating point operations per second (FLOPS), or FLOP count, as well as the memory needed to store the computed matrices. The merge process is seen as an elimination of the points on the interface of neighboring patches. For the 2-to-1 vertical merge, computing $\mathbf{S}$ requires $M^3$ FLOPS as a linear solve is necessary, and computing $\mathbf{T}$ requires $36M^3$ FLOPS. Storing $\mathbf{S}$ and $\mathbf{T}$ requires $6M^2$ and $36M^2$ numbers, respectively. For the 2-to-1 horizontal merge, computing $\mathbf{S}$ requires $8M^3$ FLOPS and computing $\mathbf{T}$ requires $128M^3$ FLOPS. Storing $\mathbf{S}$ and $\mathbf{T}$ requires $16M^2$ and $64M^2$ numbers, respectively. Thus, for a full 4-to-1 merge via two vertical merges and one horizontal merge, the total FLOP count is $210M^3$, with storage for $164M^2$ numbers. For the 4-to-1 merge, computing $\mathbf{S}$ requires $64M^3$ FLOPS and computing $\mathbf{T}$ requires $256M^3$ FLOPS. Storing $\mathbf{S}$ and $\mathbf{T}$ requires $32M^2$ and $64M^2$ numbers, respectively. Thus, for a 4-to-1 merge outlined in this paper, the total FLOP count is $320M^3$, with storage for $96M^2$ numbers.  Our 4-1 merge requires about 50\% more FLOPs  than a 2-1 merge on our finite volume mesh.  However, the 4-1 requires 70\% less storage.  We justify the higher FLOP count with the greater ease of implementation of the 4-1 merge over the 2-1, since only one type of merge algorithm is required.  

%% file: sections/implementation-details.tex
\section{Implementation Details}
\label{sec:adaptivity}

In this section, we describe the details useful in understanding how we implement our quadtree-adaptive HPS method. We discuss the patch solver implementation which wraps fast solvers. We provide details related to the mesh library \pforest which is used as a backend for the quadtree data structure.

\subsection{Patch Solvers}
\label{sub:patch_solvers}

When solving \refeq{eq:elliptic_pde} on a single leaf patch, the method used to solve the boundary value problem is called a {\em patch solver}.  In the description of the HPS method in \refsec{sec:quadtree}, the operators \Lhtau and \Linhtau are the components of the patch solver. Here, we denote this function as \texttt{PatchSolver}. The goal of the patch solver is to perform the computations in \refeq{eq:patch_leaf_solve}. The patch solver takes as input the Dirichlet data on the boundary \gext, the inhomogeneous data \ftau, and some representation of the discretization (i.e., a finite volume grid that corresponds to the patch domain). On output, \texttt{PatchSolver} returns the solution \utau of \refeq{eq:elliptic_pde} on the leaf patch.

The patch solver routine should take advantage of fast and optimized solvers for the boundary value problem \refeq{eq:elliptic_pde}. For this implementation, we wrap the FISHPACK routines \cite{swarztrauber1999fishpack} provided in the FISHPACK90 library \cite{adams2016fishpack90}. FISHPACK solves \refeq{eq:elliptic_pde} using a cyclic-reduction method, providing a fast and efficient solver for ``small'' problems.  For our cell-centered mesh, we wrap the FISHPACK routine \texttt{hwscrt} method for our \texttt{PatchSolver}.

\subsection{Building Leaf Level Operators}

At the leaf level we must construct the operator \Ttau explicitly. From \eqn{map_D2N}, we have

\begin{equation}
\Ttau = \frac{2}{h} (I - G \Lhtau). 
\end{equation}

In practice, we build \Ttau column-by-column by first solving \refeq{eq:elliptic_pde} with columns $\hat{\mathbf e}_j$ of the identity matrix $I \in \real^{4M \times 4M}$ as Dirichlet data \gtau.  Then given a solution $\mathbf u_j$ to this Dirichlet problem on the leaf, we use \eqn{vtau} to construct the corresponding column in \Ttau.

For a constant coefficient elliptic problem, we can take advantage of symmetry in \Ttau to reduce the number of calls to the \texttt{PatchSolver}. \ignore{This is due to the limited range of the Dirichlet-to-Neumann operator being a boundary operator that depends solely on the discretization of the elliptic problem.} To build \Ttau using these optimizations, we compute columns of \Ttau corresponding to one edge of the leaf patch.  The remaining columns can be filled in according to the pattern given as
\begin{align}
\textbf{T}^{\tau} =
\begin{bmatrix}
    \textbf{T}_{w,w} & -\textbf{T}_{w,e} & \textbf{T}_{w,s} & -\textbf{T}_{w,n} \\
    \textbf{T}_{w,e} & -\textbf{T}_{w,w} & \textbf{T}_{w,n} & \textbf{T}_{e,n} \\
    \textbf{T}_{w,s} & -\textbf{T}_{w,n}^T & \textbf{T}_{w,w} & -\textbf{T}_{w,e} \\
    \textbf{T}_{w,n} & \textbf{T}_{w,n}^{'} & \textbf{T}_{w,e} & -\textbf{T}_{w,w} \\
\end{bmatrix}
\end{align}
where $\textbf{T}_{w,n}^{'}$ is $\textbf{T}_{w,n}$ with reversed columns: $T_{i,j}^{'} = T_{M-i,j}$.


\subsection{Quadtrees and Adaptive Meshes}

\begin{figure}
\centering
\begin{tabular}{c c}
\smallskip
    \begin{subfigure}[t]{0.8\textwidth}
        \centering
        \includegraphics[width=\textwidth, clip=true, trim={0 150 0 150}]{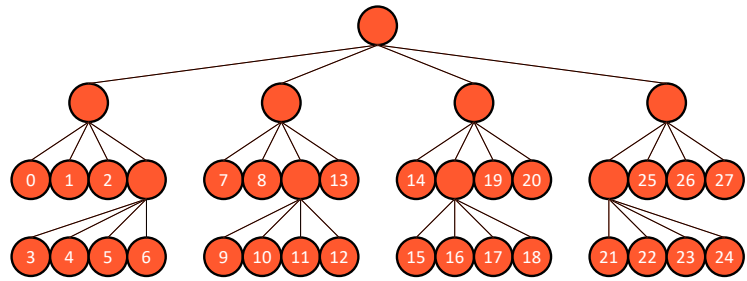}
        \caption{Leaf-level indexing of quadtree nodes}
    \end{subfigure}
    \\
    \begin{subfigure}[t]{0.8\textwidth}
        \centering
        \includegraphics[width=\textwidth, clip=true, trim={0 140 0 150}]{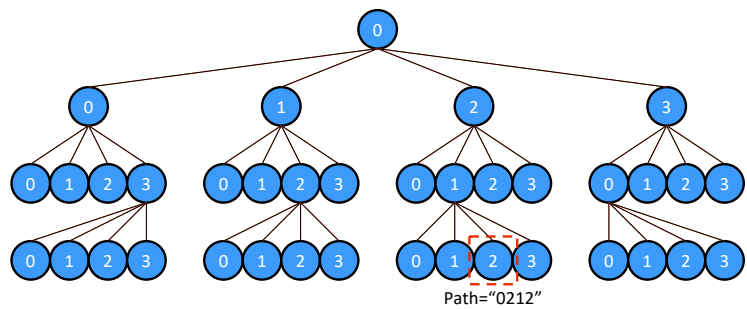}
        \caption{Path indexing of quadtree nodes}
    \end{subfigure}
\end{tabular}\\
\caption{Leaf-indexed vs. path-indexed quadtrees. In (a), only the leaves of a quadtree are indexed and stored. In (b), all nodes of the quadtree are indexed and stored according to their unique path.}
\label{fig:quadtree_indexing}
\end{figure}

The software library \pforest provides a quadtree data structure and functions to construct, store, and iterate over leaf-level quadrants in the quadtree.  However, the quadtree-adaptive HPS method requires storage for all nodes in the quadtree, including children and parent nodes in addition to leaf nodes. \ignore{This is to store the set of solution operators that act as the matrix factorization of the system matrix.} Thus, we require a new data structure that allows for data storage for all nodes in a quadtree. One approach that will later prove to be advantageous for parallel applications is a {\em path-indexed quadtree} \cite{woodward1982explicit,samet1984quadtree}. The path of a node is the unique series of directions required to traverse from the root of the tree to the node. The unique path is encoded as a string containing the sequence of nodes visited.  \reffig{fig:quadtree_indexing} illustrates the path and this unique encoding for a level 3 node

\begin{figure}
\centering
\includegraphics[width=0.4\textwidth]{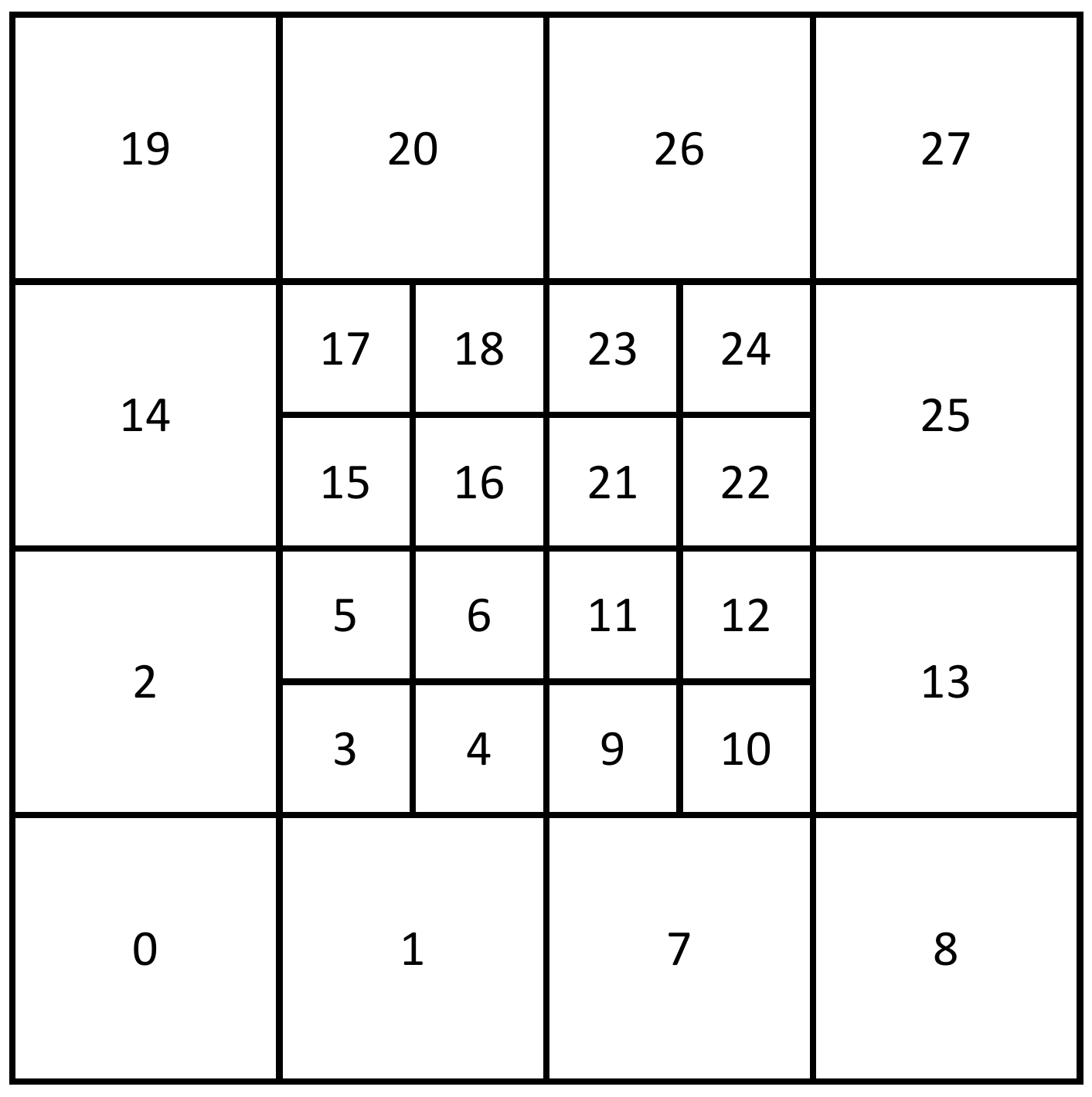}
\caption{A logically square domain is recursively broken into 4 children domains. This mesh represents refinement criteria that refines the center of the domain. The leaf level nodes are indexed according to a leaf-indexed quadtree and follow a space filling curve.}
\label{fig:adaptive_mesh}
\end{figure}

\ignore{As shown in \reffig{fig:quadtree_indexing}, each node in a family can be indexed $0$ to $3$. The highlighted node in \reffig{fig:quadtree_indexing} has a path of $0212$ which corresponds to the steps taken from the root to arrive at that node.}

\ignore{Quadtree data structures, including the path-indexed variety we use here, are not novel. One of the most notable uses of the quadtree data structure was for computer graphics (\cite{woodward1982explicit}, \cite{samet1984quadtree}). A quadtree data structure allows for hierarchical storage at different levels and is a natural application for recursive algorithms. Storing data at all nodes in a quadtree requires knowing the {\em path} of a node. The path of a node is the series of directions taken from the root of the tree to arrive at the particular node. Each node can be uniquely indexed according to the node's path. As shown in \reffig{fig:quadtree_indexing}, each node in a family can be indexed $0$ to $3$. The highlighted node in \reffig{fig:quadtree_indexing} has a path of $0212$ which corresponds to the steps taken from the root to arrive at that node.}

A path-indexed quadtree creates storage for all nodes in the tree through the use of \texttt{NodeMap}, an alias for the C++ routine \texttt{std::map<std::string, Node<T>*>}, where the template parameter \texttt{T} corresponds to a user-provided class to be stored at each node. For the quadtree-adaptive HPS method detailed in this paper, this is a class that holds patch data (grid information, matrices, and vectors). The path is computed by calling \texttt{p4est\_quadrant\_ancestor\_id}. The path-indexed quadtree wraps functions provided by \pforest to construct and iterate over nodes in the path-indexed quadtree.

Two functions provided by {\pforest} allow for a depth-first traversal of a {\pforest} quadtree: \texttt{p4est\_search\_all} and \texttt{p4est\_search\_reorder}. \texttt{p4est\_search\_all} performs a top-down search of the quadtree and provides a callback function with access to a {\pforest} quadrant data structure. This function is used to initialize the path-indexed quadtree by traversing the quadtree in a depth-first order, computing the path for each node, and allocating memory for a \texttt{Node} object in the \texttt{NodeMap}. The \texttt{p4est\-\_search\-\_all} function provides the additional capability to identify the range of processes owning (part of) a search quadrant. While the present paper covers the serial case only, an eventual parallelization will require communication up and down the tree as well as sideways to relevant processes. The function \texttt{p4est\_search\_reorder} also does a top-down search, and provides pre- and post-quadrant callback functions for accessing quadrant data.

Wrapping \texttt{p4est\_search\_reorder}, we define the {\em merge traversal} and the {\em split traversal} of a quadtree. The merge traversal iterates over a quadtree in a post-order fashion, visiting children then parents. When visiting a leaf, a leaf callback is called. When visiting parents, a family callback is used that provides access to the parent and the four children nodes. This is to ``merge'' four children nodes into a parent node, after any leaf data is assigned or computed. The split traversal iterates over a quadtree in a pre-order fashion, with callbacks to families then leaf nodes. This is to provide access to families to ``split'' one parent node into four children nodes. The leaf callback is done last in the split traversal and is used to compute leaf level data once the entire tree has been traversed. The algorithm for the callback function passed to \texttt{p4est\_search\_reorder} is provided in \refalg{alg:quadtree_callback}.

\begin{algorithm}
\caption{\texttt{QuadtreeCallback} Function}
\begin{algorithmic}[0]
    \Require Functions \texttt{LeafCallback}(\texttt{leaf\_node}) \& \texttt{FamilyCallback}(\texttt{parent\_node, children\_nodes})
    \State Compute \texttt{path} from \texttt{p4est\_quadrant\_ancestor\_id}
    \State Let \texttt{node = node\_map[path]}
    \If{\texttt{node} is a leaf} \Comment{Node is a leaf, call leaf callback}
        \State Call \texttt{LeafCallback(node)}
    \Else \Comment{Node is a branch, get children and call family callback}
        \For{i = 0, 1, 2, 3}
            \State \texttt{children\_nodes[i] = map[path + string(i)]}
        \EndFor
        \State Call \texttt{FamilyCallback(node, children\_nodes)}
    \EndIf
\end{algorithmic}
\label{alg:quadtree_callback}
\end{algorithm}

For the quadtree-adaptive HPS method, the leaf callback for the merge traversal is \texttt{BuildD2N}, which solves \refeq{eq:map_D2N} and computes a leaf level Dirichlet-to-Neumann matrix. The family callback wraps the algorithms found in \refsec{sub:4-to-1merge}, which performs the 4-to-1 merge.

The family callback for the split traversal wraps \refeq{eq:solve_eqn}, which applies the solution operator to map parent Dirichlet data to children Dirichlet data (the 1-to-4 split). The leaf callback wraps \texttt{PatchSolver} to solve \refeq{eq:elliptic_pde} on a leaf level using fast solution methods.

%% file: sections/results.tex
\section{Numerical Results}
\label{sec:results}

To test and verify the implementation of the quadtree-adaptive HPS method, we solve two Poisson equations and one Helmholtz equation. For each, we present error, timing, and memory usage results and discuss the performance of our implemented method.

For all of our examples, we discretize the Laplace operator at the patch level using a second-order, 5-point stencil on a finite volume (cell-centered) mesh. The code, along with numerical experiments below, is stored in the GitHub repository EllipticForest \cite{chipman2023elliptic}. EllipticForest is written primarily in C++, with wrappers to FORTRAN routines to call FISHPACK and LAPACK (\cite{anderson1999lapack}) for any dense linear algebra operations. The mesh and solution are output into an unstructured VTK mesh file \cite{vtkBook} and are visualized with the VisIt software \cite{HPV:VisIt}. All tests were run on a 2021 MacBook Pro with an M1 Pro CPU and 32 GB of RAM.

\subsection{Poisson Equation 1}
\label{sub:example_one}

We solve the following boundary value problem:
\begin{align}
    \nabla^2 u(x,y) = -(\sin(x) + \sin(y))
\end{align}
on the square domain $\Omega = [-10, 10] \times [-10, 10]$, subject to Dirichlet boundary conditions $u(x,y) = g(x,y)$ on the boundary which is computed according to the exact solution
\begin{align}
    u_{exact}(x,y) &= \sin(x) + \sin(y).
\end{align}

\begin{figure}
    \centering
    \includegraphics[width=0.75\textwidth, trim={0 100 0 0}]{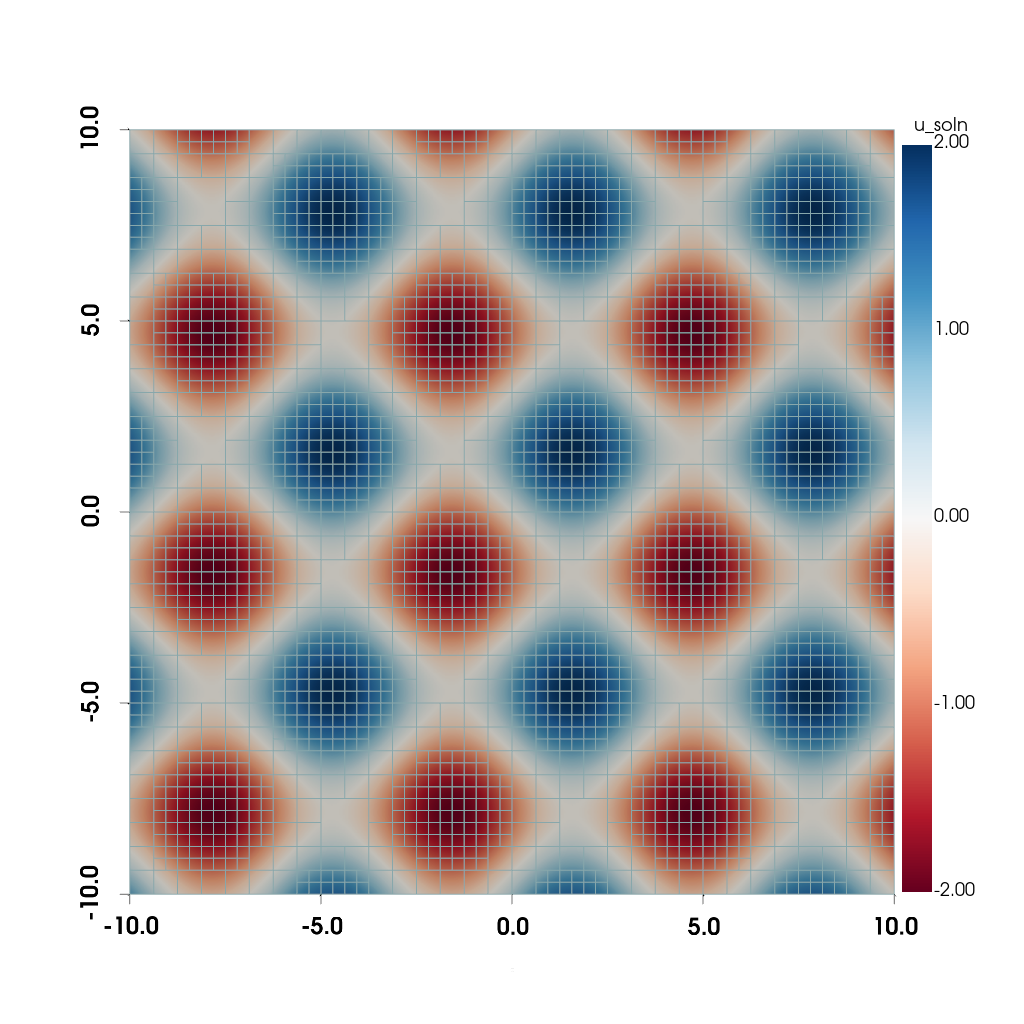}
    \caption{The computed solution and mesh for the Poisson problem \refsec{sub:example_one}.  Patch size for this plot is $16 \times 16$ and mesh is refined to level 7.  Refinement criteria is based on the magnitude of the right hand side function $f(x,y)$.}
    \label{fig:poisson_plot}
\end{figure}

For the refinement criteria, we refine according to the right-hand side function $f(x,y)$, which corresponds to the curvature of the solution. We set a refinement threshold of $1.2$ and refine a patch when $f(x,y) > 1.2$ for any $x,y$ in a patch. This results in a mesh and solution that can be found in \reffig{fig:poisson_plot}.

{\bf Results and Discussion}
Tables \ref{tab:poisson_error} and \ref{tab:poisson_timing} show the error, timing, and memory results for the current implementation on this Poisson equation. \reftab{tab:poisson_error} shows results for both a uniformly refined mesh (a mesh without any coarse-fine interfaces or local adaptivity) and results for the adaptive mesh case. For the uniform case, we get the expected second order convergence in both the $L_{\infty}$ and $L_1$ norms. The adaptive case mostly shows second order convergence, except for a few cases where the refinement between successive levels results in a smaller jump in error than expected. \reftab{tab:poisson_timing} shows timing and memory results for the same case. Here, the difference between the uniform and adaptive case is highlighted. The adaptive case gives a $4.5$ times speed up for the build stage, and nearly a $20$ times speed up for the solve stage. Memory used to store the quadtree and operators is also significantly reduced with the adaptive case.

\begin{table}
    \caption{Convergence analysis for Poisson's equation. The upper part shows convergence for a uniformly refined mesh, while the lower part shows convergence for an adaptively refined mesh. $M$ is the size of the grid on each leaf patch, $L_{\text{max}}$ is the maximum level of refinement, $R_{\text{eff}}$ is the effective resolution for a uniformly refined mesh, DOFs is the total degrees of freedom (i.e., total mesh points), $L_{\infty}$ error is the infinity norm error, $L_{\infty}$ order is the infinity norm convergence order, $L_1$ error is the $1^{\text{st}}$ norm error, and $L_1$ order is the $1^{\text{st}}$ norm convergence order.}
    \centering
    \sisetup{
        table-alignment-mode=format
    }
    \begin{tabular}{
        |
        S   
        S[table-column-width=0.8cm]   
        S[table-text-alignment=right, table-column-width=0.7cm]   
        S[table-text-alignment=right, table-column-width=1.4cm]   
        S[scientific-notation=true, round-mode=places, round-precision=2, table-column-width=1.65cm]   
        S[scientific-notation=false, exponent-mode=fixed, round-mode=places, round-precision=2, table-column-width=1.5cm]   
        S[scientific-notation=true, round-mode=places, round-precision=2]   
        S[scientific-notation=false, exponent-mode=fixed, round-mode=places, round-precision=2]   
        |
    }
\hline
{M} & {$L_{\text{max}}$} & {$R_{\text{eff}}$} & {DOFs} & {$L_{\infty}$ Error} & {$L_{\infty}$ Order} & {$L_1$ Error} & {$L_1$ Order} \\
\hline
\num{16} & \num{4} & \num{256} & \num{65536} & \num{1.1146466297640E-03} & \num{2.0002722123237E+00} & \num{3.5892082197986E-04} & \num{2.0009066196875E+00} \\
\num{16} & \num{5} & \num{512} & \num{262144} & \num{2.7855514231123E-04} & \num{2.0005515585727E+00} & \num{8.9717104596707E-05} & \num{2.0002106534890E+00} \\
\num{16} & \num{6} & \num{1024} & \num{1048576} & \num{6.9642847198903E-05} & \num{1.9999158585307E+00} & \num{2.2428541266291E-05} & \num{2.0000472698832E+00} \\
\num{16} & \num{7} & \num{2048} & \num{4194304} & \num{1.7410906046011E-05} & \num{1.9999839043577E+00} & \num{5.6070914295010E-06} & \num{2.0000112920243E+00} \\
\num{32} & \num{4} & \num{512} & \num{262144} & \num{2.7855514234232E-04} & \num{2.0005515584215E+00} & \num{8.9717104608167E-05} & \num{2.0002106533137E+00} \\
\num{32} & \num{5} & \num{1024} & \num{1048576} & \num{6.9642847325024E-05} & \num{1.9999158560790E+00} & \num{2.2428541313515E-05} & \num{2.0000472670298E+00} \\
\num{32} & \num{6} & \num{2048} & \num{4194304} & \num{1.7410906499649E-05} & \num{1.9999838693813E+00} & \num{5.6070916187308E-06} & \num{2.0000112463735E+00} \\
\hline
\num{16} & \num{4} & \num{256} & \num{64000} & \num{2.5910957992683E-03} & \num{7.8329626881446E-01} & \num{4.4586170147130E-04} & \num{1.6879759592817E+00} \\
\num{16} & \num{5} & \num{512} & \num{194560} & \num{6.6263501289099E-04} & \num{1.9672760156010E+00} & \num{1.2477150864153E-04} & \num{1.8373077457909E+00} \\
\num{16} & \num{6} & \num{1024} & \num{569344} & \num{6.8094879278413E-04} & \num{-3.9331876110758E-02} & \num{1.0216613949030E-04} & \num{2.8837140596648E-01} \\
\num{16} & \num{7} & \num{2048} & \num{1984000} & \num{1.7141451745961E-04} & \num{1.9900570123994E+00} & \num{3.7577025980665E-05} & \num{1.4429943340785E+00} \\
\num{32} & \num{4} & \num{512} & \num{256000} & \num{6.6189922864102E-04} & \num{7.5190291839197E-01} & \num{1.1198962819648E-04} & \num{1.6803004955197E+00} \\
\num{32} & \num{5} & \num{1024} & \num{784384} & \num{1.6733508690897E-04} & \num{1.9838716073808E+00} & \num{3.1111362539336E-05} & \num{1.8478516397654E+00} \\
\num{32} & \num{6} & \num{2048} & \num{2289664} & \num{1.7139947044398E-04} & \num{-3.4622669955039E-02} & \num{2.5679498208238E-05} & \num{2.7682456821605E-01} \\
\hline
    \end{tabular}
    \label{tab:poisson_error}
\end{table}

\begin{table}
    \caption{Timing and memory results for Poisson's equation. The first four rows show results for the uniformly refined mesh, while the second set of rows shows results for the adaptively refined mesh. The results here are for a patch size of $16 \times 16$. $L_{\mbox{max}}$ is the maximum level of refinement, $R_{\text{eff}}$ is the effective resolution, DOFs is the total degrees of freedom, $T_{\text{build}}$ is the time in seconds for the build stage, $T_{\text{upwards}}$ is the time in seconds for the upwards stage, $T_{\text{solve}}$ is the time in seconds for the solve stage, and $S$ is the memory storage in megabytes to store the quadtree and all data matrices stored in each node of the quadtree.}
    \centering
    \sisetup{
        table-alignment-mode=format
    }
    \begin{tabular}{
        |
        S   
        S[table-text-alignment=right]   
        S[table-text-alignment=right]   
        S[table-text-alignment=right, scientific-notation=false, round-mode=places, round-precision=3]   
        S[table-text-alignment=right, scientific-notation=false, exponent-mode=fixed, round-mode=places, round-precision=3]   
        S[table-text-alignment=right, scientific-notation=false, round-mode=places, round-precision=3]   
        S[table-text-alignment=right, scientific-notation=false, round-mode=places, round-precision=1]   
        |
    }
\hline
{$L_{\text{max}}$} & {$R_{\text{eff}}$} & {DOFs} & {$T_{\text{build}}$ (sec)} & {$T_{\text{upwards}}$ (sec)} & {$T_{\text{solve}}$ (sec)} & {$S$ (MB)} \\
\hline
\num{4} & \num{256} & \num{65536} & \num{2.031286957999990} & \num{3.66539416999999E-01} & \num{0.211029084000000} & \num{81.548497200012200} \\
\num{5} & \num{512} & \num{262144} & \num{8.525974832999990} & \num{1.77923349999999E+00} & \num{1.073760000000000} & \num{398.233067512512000} \\
\num{6} & \num{1024} & \num{1048576} & \num{39.353658666999900} & \num{9.00205208400000E+00} & \num{4.258692584000000} & \num{1881.010411262510000} \\
\num{7} & \num{2048} & \num{4194304} & \num{172.189544875000000} & \num{4.20349420419999E+01} & \num{20.559359582999900} & \num{8676.197911262510000} \\
\hline
\num{4} & \num{256} & \num{64000} & \num{1.502292208999990} & \num{3.93616458000000E-01} & \num{0.098160250000000} & \num{52.977078437805100} \\
\num{5} & \num{512} & \num{194560} & \num{3.718632625000000} & \num{9.24946582999999E-01} & \num{0.290812499999999} & \num{140.070134162902000} \\
\num{6} & \num{1024} & \num{569344} & \num{10.594630750000000} & \num{2.22070095800000E+00} & \num{0.478022832999999} & \num{368.462132453918000} \\
\num{7} & \num{2048} & \num{1984000} & \num{38.976778125000000} & \num{8.09454862499999E+00} & \num{1.838669208000000} & \num{1443.861859321590000} \\
\hline
    \end{tabular}
    \label{tab:poisson_timing}
\end{table}

\subsection{Poisson Equation 2 (Polar-Star Problem)}

\begin{figure}
    \centering
    \includegraphics[width=0.75\textwidth, trim={0 100 0 0}]{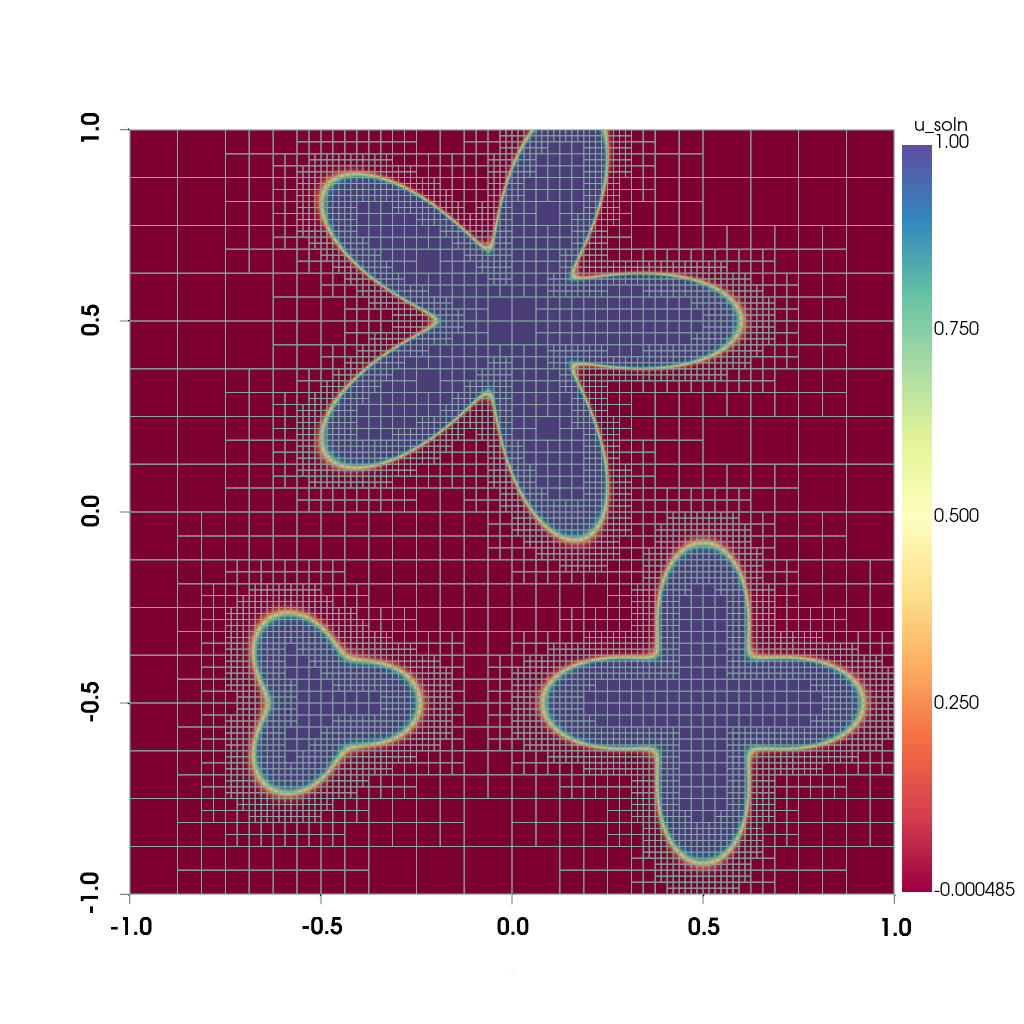}
    \caption{The computed solution and mesh for the Polar Star Poisson problem. Each patch has a $16 \times 16$ cell-centered grid. The mesh is refined according to the right-hand side and is refined with 8 levels of refinement.}
    \label{fig:polar_star_plot}
\end{figure}

For this second test, we solve a ``polar star'' Poisson problem. The problem is created via a method of manufactured solutions and is engineered to have highly local curvature from the load function. The resulting solution is a collection of polar stars with user specified number of points and radii of curvature. This problem highlights the use of an adaptive mesh to solve the elliptic equation. The exact solution we attempt to reconstruct is the following:
\begin{align}
    u(x,y) = \frac{1}{2} \sum_{i=1}^{N} 1 - \tanh \left(\frac{r(x,y)-r_{0,i} \left(r_{1,i} \cos \left(n \theta(x,y)\right)+1\right)}{\epsilon }\right)
\end{align}
Computing the Laplacian analytically yields the right-hand side to Poisson's equation. Thus, the polar star Poisson problem is defined as follows:
\begin{align}
    \nabla^2 u(x,y) = \sum_{i=1}^N -\frac{s_{1,i}(x,y) + s_{2,i}(x, y)}{r(x,y)^2} - s_{3,i}(x,y) + s_{4,i}(x,y)
\end{align}
with
\begin{align*}
    s_{1,i}(x,y) &= \frac{p(x,y)^2 \tanh \left(\phi(x,y)\right) \text{sech}^2\left(\phi(x,y)\right)}{\epsilon ^2} \\
    s_{2,i}(x,y) &= -\frac{n^2 r_{0,i} r_{1,i} \cos (n \theta(x,y)) \text{sech}^2\left(\phi(x,y)\right)}{2 \epsilon } \\
    s_{3,i}(x,y) &= \frac{\tanh \left(\phi(x,y)\right) \text{sech}^2\left(\phi(x,y)\right)}{\epsilon ^2} \\
    s_{4,i}(x,y) &= \frac{\text{sech}^2\left(\phi(x,y)\right)}{2 r(x,y) \epsilon } \\
    p(x,y) &= n r_{0,i} r_{1,i} \sin (n \theta(x,y)) \\
    \phi(x,y) &= \frac{r(x,y)-r_{0,i} (r_{1,i} \cos (n \theta(x,y))+1)}{\epsilon}
\end{align*}
and where $i=1, ..., N_{polar}$ and $N_{polar}$ is the number of polar stars. Each polar star has a center $(x_0, y_0)$, inner and outer radii $r_0, r_1$, and the number of arms per polar star $n$. The radius and angle have the standard polar transforms:
\begin{align}
    r(x,y) &= \sqrt{(x - x_0)^2 + (y - y_0)^2} \\
    \theta(x,y) &= \tan^{-1}\Big(\frac{y - y_0}{x - x_0}\Big)
\end{align}
\reftab{table:polar_star_parameters} has the parameters used in this case study. \reffig{fig:polar_star_plot} shows the resulting mesh and solution. \reftab{tab:polar_star_results} shows the error, timing, and memory results for the polar star Poisson problem.
\begin{table}[ht]
    \begin{center}
        \caption{Polar Star Poisson Problem Parameters}
        \begin{tabular}{|c|c|c|c|c|c|}
            \hline
            $i$ & $x_0$ & $y_0$ & $r_0$ & $r_1$ & $n$ \\
            \hline
            $1$ & $-0.5$ & $-0.5$ & $0.2$ & $0.3$ & $3$ \\
            $2$ & $0.5$ & $-0.5$ & $0.3$ & $0.4$ & $4$ \\
            $3$ & $0$ & $0.5$ & $0.4$ & $0.5$ & $5$ \\
            \hline
        \end{tabular}
        \label{table:polar_star_parameters}
    \end{center}
\end{table}

\begin{table}
    \caption{Convergence analysis for the Polar Star Poisson problem. The upper part shows convergence for a uniformly refined mesh, while the lower part shows convergence for an adaptively refined mesh. $M$ is the size of the grid on each leaf patch, $L_{\text{max}}$ is the maximum level of refinement, $R_{\text{eff}}$ is the effective resolution for a uniformly refined mesh, DOFs is the total degrees of freedom (i.e., total mesh points), $L_{\infty}$ error is the infinity norm error, $L_{\infty}$ order is the infinity norm convergence order, $L_1$ error is the $1^{\text{st}}$ norm error, and $L_1$ order is the $1^{\text{st}}$ norm convergence order.}
    \centering
    \sisetup{
        table-alignment-mode=format
    }
    \begin{tabular}{
        |
        S   
        S[table-column-width=0.8cm]   
        S[table-text-alignment=right, table-column-width=0.7cm]   
        S[table-text-alignment=right, table-column-width=1.4cm]   
        S[scientific-notation=true, round-mode=places, round-precision=2, table-column-width=1.65cm]   
        S[scientific-notation=false, exponent-mode=fixed, round-mode=places, round-precision=2, table-column-width=1.5cm]   
        S[scientific-notation=true, round-mode=places, round-precision=2]   
        S[scientific-notation=false, exponent-mode=fixed, round-mode=places, round-precision=2]   
        |
    }
\hline
{M} & {$L_{\text{max}}$} & {$R_{\text{eff}}$} & {DOFs} & {$L_{\infty}$ Error} & {$L_{\infty}$ Order} & {$L_1$ Error} & {$L_1$ Order} \\
\hline
\num{16} & \num{4} & \num{256} & \num{65536} & \num{1.56199870806249E+00} & \num{3.03980331387296E+00} & \num{1.65069071070283E-01} & \num{4.25135175327624E+00} \\
\num{16} & \num{5} & \num{512} & \num{262144} & \num{2.79704505790678E-02} & \num{5.80334595545976E+00} & \num{9.84718059820731E-04} & \num{7.38914339546720E+00} \\
\num{16} & \num{6} & \num{1024} & \num{1048576} & \num{5.25422685067089E-03} & \num{2.41235309940473E+00} & \num{7.39997994016763E-05} & \num{3.73411745253242E+00} \\
\num{16} & \num{7} & \num{2048} & \num{4194304} & \num{1.28369442835374E-03} & \num{2.03317666696840E+00} & \num{1.84035492646184E-05} & \num{2.00753733205534E+00} \\
\num{32} & \num{3} & \num{256} & \num{65536} & \num{1.56199870806260E+00} & \num{3.03980331387287E+00} & \num{1.65069071070315E-01} & \num{4.25135175327601E+00} \\
\num{32} & \num{4} & \num{512} & \num{262144} & \num{2.79704505792760E-02} & \num{5.80334595544912E+00} & \num{9.84718059909330E-04} & \num{7.38914339533767E+00} \\
\num{32} & \num{5} & \num{1024} & \num{1048576} & \num{5.25422685134335E-03} & \num{2.41235309923083E+00} & \num{7.39997994074657E-05} & \num{3.73411745254936E+00} \\
\num{32} & \num{6} & \num{2048} & \num{4194304} & \num{1.28369443102049E-03} & \num{2.03317666415598E+00} & \num{1.84035493232030E-05} & \num{2.00753732757563E+00} \\
\hline
\num{16} & \num{4} & \num{256} & \num{50944} & \num{1.56155552227466E+00} & \num{3.04021270772076E+00} & \num{1.64990394749018E-01} & \num{4.25203954412867E+00} \\
\num{16} & \num{5} & \num{512} & \num{171520} & \num{7.12119724134531E-02} & \num{4.45472024347345E+00} & \num{1.03166637134095E-03} & \num{7.32126173153682E+00} \\
\num{16} & \num{6} & \num{1024} & \num{476416} & \num{9.14563250417531E-02} & \num{-3.60963136314842E-01} & \num{1.78815639727227E-04} & \num{2.52843166634489E+00} \\
\num{16} & \num{7} & \num{2048} & \num{1587712} & \num{2.81786611011458E-02} & \num{1.69847988484306E+00} & \num{6.32775510074655E-05} & \num{1.49870725401785E+00} \\
\num{32} & \num{3} & \num{256} & \num{65536} & \num{1.56199870806260E+00} & \num{3.03980331387287E+00} & \num{1.65069071070315E-01} & \num{4.25135175327601E+00} \\
\num{32} & \num{4} & \num{512} & \num{203776} & \num{2.79706867876203E-02} & \num{5.80333377204982E+00} & \num{9.84314973348746E-04} & \num{7.38973407206102E+00} \\
\num{32} & \num{5} & \num{1024} & \num{701440} & \num{2.33901397925928E-02} & \num{2.58015193873839E-01} & \num{8.15186205620300E-05} & \num{3.59391849705353E+00} \\
\num{32} & \num{6} & \num{2048} & \num{1921024} & \num{2.81786212280739E-02} & \num{-2.68700538645513E-01} & \num{3.10231851870699E-05} & \num{1.39378282159843E+00} \\
\hline
    \end{tabular}
    \label{tab:polar_star_results}
\end{table}

\begin{table}
    \caption{Timing and memory results for the Polar Star Poisson problem. The first eight rows show results for the uniformly refined mesh, while the second set of rows shows results for the adaptively refined mesh. The results here are for a patch size of $16 \times 16$. $L_{\text{max}}$ is the maximum level of refinement, $R_{\text{eff}}$ is the effective resolution, DOFs is the total degrees of freedom, $T_{\text{build}}$ is the time in seconds for the build stage, $T_{\text{upwards}}$ is the time in seconds for the upwards stage, $T_{\text{solve}}$ is the time in seconds for the solve stage, and $S$ is the memory storage in megabytes to store the quadtree and all data matrices stored in each node of the quadtree.}
    \centering
    \sisetup{
        table-alignment-mode=format
    }
    \begin{tabular}{
        |
        S[table-column-width=0.8cm]   
        S[table-text-alignment=right, table-column-width=0.7cm]   
        S[table-text-alignment=right, table-column-width=1.4cm]   
        S[table-text-alignment=right, scientific-notation=false, exponent-mode=fixed, round-mode=places, round-precision=2]   
        S[table-text-alignment=right, scientific-notation=false, exponent-mode=fixed, round-mode=places, round-precision=2]   
        S[table-text-alignment=right, scientific-notation=false, exponent-mode=fixed, round-mode=places, round-precision=2]   
        S[table-text-alignment=right, scientific-notation=false, exponent-mode=fixed, round-mode=places, round-precision=2]   
        |
    }
\hline
{$L_{\text{max}}$} & {$R_{\text{eff}}$} & {DOFs} & {$T_{\text{build}}$ (sec)} & {$T_{\text{upwards}}$ (sec)} & {$T_{\text{solve}}$ (sec)} & {$S$ (MB)} \\
\hline
\num{4} & \num{256} & \num{65536} & \num{1.60620216600000E+00} & \num{4.159940420000E-01} & \num{2.08710665999999E-01} & \num{8.15484972000122E+01} \\
\num{5} & \num{512} & \num{262144} & \num{7.76756095899999E+00} & \num{1.969481791000E+00} & \num{9.65551166999999E-01} & \num{3.98233067512512E+02} \\
\num{6} & \num{1024} & \num{1048576} & \num{3.56706281669999E+01} & \num{9.239545625000E+00} & \num{4.09653733299999E+00} & \num{1.88101041126251E+03} \\
\num{7} & \num{2048} & \num{4194304} & \num{1.65943133041000E+02} & \num{4.246258775000E+01} & \num{1.93515662920000E+01} & \num{8.67619791126251E+03} \\
\hline
\num{4} & \num{256} & \num{50944} & \num{9.29590041999999E-01} & \num{2.385920420000E-01} & \num{8.11892919999999E-02} & \num{3.35243158340454E+01} \\
\num{5} & \num{512} & \num{171520} & \num{3.32502091600000E+00} & \num{7.851652910000E-01} & \num{1.44476500000000E-01} & \num{1.15864575386047E+02} \\
\num{6} & \num{1024} & \num{476416} & \num{9.06574687500000E+00} & \num{1.937123083000E+00} & \num{3.00085959000000E-01} & \num{2.95106141090393E+02} \\
\num{7} & \num{2048} & \num{1587712} & \num{3.07760261670000E+01} & \num{7.226564708000E+00} & \num{1.12108816699999E+00} & \num{1.07681472492218E+03} \\
\hline
    \end{tabular}
    \label{tab:polar_star_timing}
\end{table}

{\bf Results and Discussion}
The polar star Poisson problem highlights the benefits of an adaptive mesh. As shown in \reffig{fig:polar_star_plot}, there is room for significant speed up on the areas outside each polar star. The localized curvature is sufficiently captured by this adaptive scheme. \reftab{tab:polar_star_results} show the error analysis for this problem. Because of the large curvature at the boundaries of a polar star and our use of a second-order accurate solver, the $L_{\infty}$ error is larger than the first Poisson equation from \refsec{sub:example_one}, but we can still achieve 6 digits of accuracy in the $L_1$ norm. The timing results shown in \reftab{tab:polar_star_timing} show significant speed up between the uniform and adaptive implementations. The build time has a $5.5$ times speed up, the upwards stage has a $5.8$ times speed up, and the solve stage has a $17$ times speed up. Plus, the memory for the adaptive case is about $1/8th$ of the memory required for the uniform case.

\subsection{Helmholtz's Equation}
\label{sub:helmholtz_equation}

For our third test, we solve the boundary value problem
\begin{align}
    \nabla^2 u(x,y) + \lambda u(x,y) = f(x,y),
\end{align}
on the square domain $\Omega = [-0.5, 0.5] \times [-0.5, 0.5]$ and subject to Dirichlet boundary conditions supplied via an exact solution discussed below.

\begin{figure}
    \centering
    \begin{tabular}{c c}
    \smallskip
        \begin{subfigure}[t]{0.455\textwidth}
            \centering
            \includegraphics[width=1.00\textwidth, clip=true, trim={0 80 0 0}]{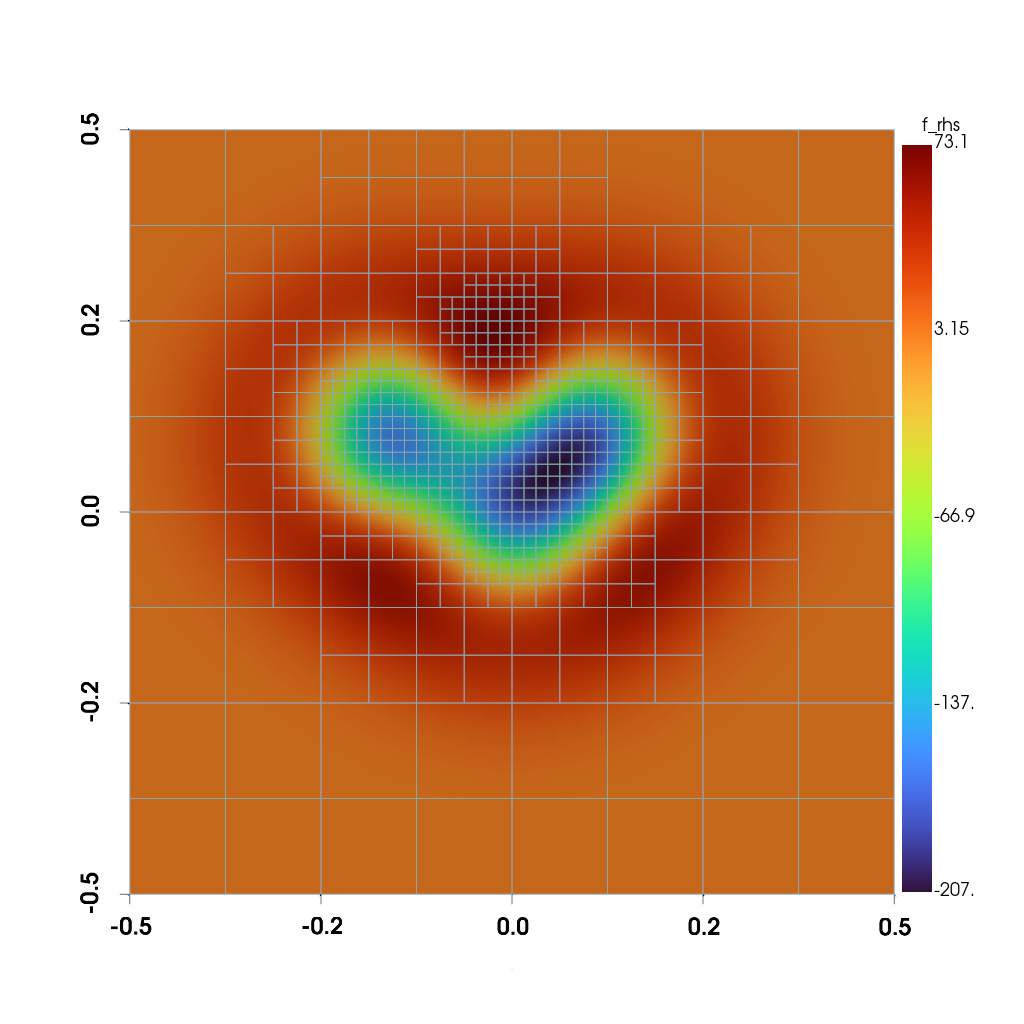}
        \end{subfigure}
        &
        \begin{subfigure}[t]{0.455\textwidth}
            \centering
            \includegraphics[width=1.00\textwidth, clip=true, trim={0 80 0 0}]{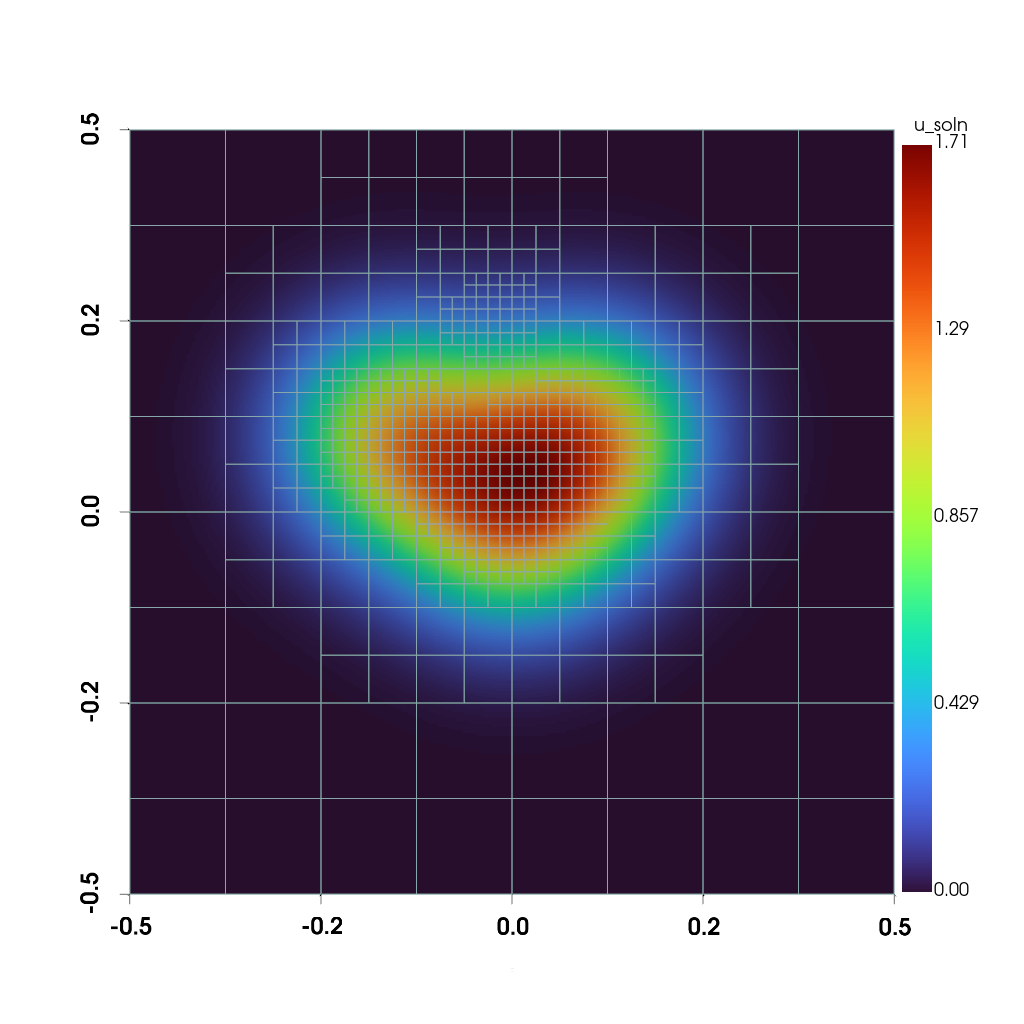}
        \end{subfigure}
    \end{tabular}\\
    \caption{Mesh and plot of the Helmholtz problem in \refsec{sub:helmholtz_equation}. The left plot is the right-hand side and the right plot is the computed solution.}
    \label{fig:helmholtz_plots}
\end{figure}

To provide an analytical solution to compare against, we use a problem provided in \cite{cheng2006adaptive}, and given by
\begin{align}
    u(\textbf{x}) = \sum_{i=1}^{3} e^{-\alpha |\textbf{x} - \textbf{x}_i|^2}
\end{align}
with $\lambda=0.01$, $\alpha=50$, $\textbf{x}_1 = (0.1, 0.1)$, $\textbf{x}_2 = (0, 0)$, and $\textbf{x}_3 = (-0.15, 0.1)$. The right-hand side is computed analytically using Mathematica \cite{Mathematica}. We set the threshold for refinement to 60. The solution and right-hand side function for $16 \times 16$ patches are plotted in \reffig{fig:helmholtz_plots}.

{\bf Results and Discussion}
As demonstrated with this example, the implemented HPS solves Helmholtz equations as well with the expected second order accuracy. See \reftab{tab:helmholtz_results} for error and convergence analysis. The adaptive mesh is able to successfully capture the curvature and reduce the work required to solve this equation. \reftab{tab:helmholtz_timing} show timing and memory usage results, which show significant speedup due to good mesh adaptation. We achieve a $17$ times speedup for the build and upwards stage and a $57$ times speedup for the solve stage. Memory is also impressive for this particular problem, with solving the same problem with similar error with $4.6\%$ of the required memory.

\begin{table}
    \caption{Convergence analysis for the Helmholtz problem. The first eight rows show convergence for a uniformly refined mesh, while the bottom eight rows show convergence for an adaptively refined mesh. $M$ is the size of the grid on each leaf patch, $L_{\text{max}}$ is the maximum level of refinement, $R_{\text{eff}}$ is the effective resolution for a uniformly refined mesh, DOFs is the total degrees of freedom (i.e., total mesh points), $L_{\infty}$ error is the infinity norm error, $L_{\infty}$ order is the infinity norm convergence order, $L_1$ error is the $1^{\text{st}}$ norm error, and $L_1$ order is the $1^{\text{st}}$ norm convergence order.}
    \centering
    \sisetup{
        table-alignment-mode=format
    }
    \begin{tabular}{
        |
        S   
        S[table-column-width=0.8cm]   
        S[table-text-alignment=right, table-column-width=0.7cm]   
        S[table-text-alignment=right, table-column-width=1.4cm]   
        S[scientific-notation=true, round-mode=places, round-precision=2, table-column-width=1.65cm]   
        S[scientific-notation=false, exponent-mode=fixed, round-mode=places, round-precision=2, table-column-width=1.5cm]   
        S[scientific-notation=true, round-mode=places, round-precision=2]   
        S[scientific-notation=false, exponent-mode=fixed, round-mode=places, round-precision=2]   
        |
    }
\hline
{M} & {$L_{\text{max}}$} & {$R_{\text{eff}}$} & {DOFs} & {$L_{\infty}$ Error} & {$L_{\infty}$ Order} & {$L_1$ Error} & {$L_1$ Order} \\
\hline
\num{16} & \num{3} & \num{128} & \num{16384} & \num{8.11856103331010E-04} & \num{1.99880860589163E+00} & \num{8.79086298739653E-05} & \num{2.00167127817603E+00} \\
\num{16} & \num{4} & \num{256} & \num{65536} & \num{2.02984786720650E-04} & \num{1.99985243641999E+00} & \num{2.19695214703194E-05} & \num{2.00050135378016E+00} \\
\num{16} & \num{5} & \num{512} & \num{262144} & \num{5.07454266427398E-05} & \num{2.00002189203670E+00} & \num{5.49192113777069E-06} & \num{2.00012063189787E+00} \\
\num{16} & \num{6} & \num{1024} & \num{1048576} & \num{1.26864921798919E-05} & \num{1.99998458881074E+00} & \num{1.37295318659519E-06} & \num{2.00002847405489E+00} \\
\num{16} & \num{7} & \num{2048} & \num{4194304} & \num{3.17161259033582E-06} & \num{2.00000475557083E+00} & \num{3.43237939677034E-07} & \num{2.00000150042016E+00} \\
\num{32} & \num{3} & \num{256} & \num{65536} & \num{2.02984786837223E-04} & \num{1.99985243560645E+00} & \num{2.19695214601884E-05} & \num{2.00050135443361E+00} \\
\num{32} & \num{4} & \num{512} & \num{262144} & \num{5.07454269149665E-05} & \num{2.00002188512581E+00} & \num{5.49192111200352E-06} & \num{2.00012063800147E+00} \\
\num{32} & \num{5} & \num{1024} & \num{1048576} & \num{1.26864947160854E-05} & \num{1.99998430813683E+00} & \num{1.37295295907976E-06} & \num{2.00002870635855E+00} \\
\num{32} & \num{6} & \num{2048} & \num{4194304} & \num{3.17160773866120E-06} & \num{2.00000725090321E+00} & \num{3.43238375314458E-07} & \num{1.99999943028095E+00} \\
\hline
\num{16} & \num{3} & \num{128} & \num{8704} & \num{1.35786298229745E-03} & \num{1.25676664286521E+00} & \num{1.62110560935590E-04} & \num{1.11876990262527E+00} \\
\num{16} & \num{4} & \num{256} & \num{22528} & \num{1.64905762929112E-03} & \num{-2.80303908147837E-01} & \num{8.45781698947318E-05} & \num{9.38620831715358E-01} \\
\num{16} & \num{5} & \num{512} & \num{54784} & \num{2.38573635354355E-03} & \num{-5.32792803144748E-01} & \num{3.38667165953099E-05} & \num{1.32041722042279E+00} \\
\num{16} & \num{6} & \num{1024} & \num{163072} & \num{5.94628321314072E-04} & \num{2.00437453684095E+00} & \num{1.31283531083692E-05} & \num{1.36718217508683E+00} \\
\num{16} & \num{7} & \num{2048} & \num{485632} & \num{7.32427776781507E-04} & \num{-3.00698326888675E-01} & \num{1.91363184831188E-05} & \num{-5.43627357312482E-01} \\
\num{32} & \num{3} & \num{256} & \num{34816} & \num{3.39042440507864E-04} & \num{1.25975816324197E+00} & \num{4.06849021440135E-05} & \num{1.11151127906479E+00} \\
\num{32} & \num{4} & \num{512} & \num{90112} & \num{4.02109898754332E-04} & \num{-2.46123973731869E-01} & \num{2.13036585418219E-05} & \num{9.33392310591107E-01} \\
\num{32} & \num{5} & \num{1024} & \num{222208} & \num{5.85657916098325E-04} & \num{-5.42468378289773E-01} & \num{8.31054203108411E-06} & \num{1.35808672932660E+00} \\
\num{32} & \num{6} & \num{2048} & \num{664576} & \num{1.50576960741055E-04} & \num{1.95955718461569E+00} & \num{3.71801364039699E-06} & \num{1.16041051254945E+00} \\
\hline
    \end{tabular}
    \label{tab:helmholtz_results}
\end{table}

\begin{table}
    \caption{Timing and memory results for the Helmholtz problem. The upper part shows results for the uniformly refined mesh, while the lower part shows results for the adaptively refined mesh. The results here are for a patch size of $16 \times 16$. $L_{\text{max}}$ is the maximum level of refinement, $R_{\text{eff}}$ is the effective resolution, DOFs is the total degrees of freedom, $T_{\text{build}}$ is the time in seconds for the build stage, $T_{\text{upwards}}$ is the time in seconds for the upwards stage, $T_{\text{solve}}$ is the time in seconds for the solve stage, and $S$ is the memory storage in megabytes to store the quadtree and all data matrices stored in each node of the quadtree.}
    \centering
    \sisetup{
        table-alignment-mode=format
    }
    \begin{tabular}{
        |
        S[table-column-width=0.8cm]   
        S[table-text-alignment=right, table-column-width=0.7cm]   
        S[table-text-alignment=right, table-column-width=1.4cm]   
        S[table-text-alignment=right, scientific-notation=false, exponent-mode=fixed, round-mode=places, round-precision=2]   
        S[table-text-alignment=right, scientific-notation=false, exponent-mode=fixed, round-mode=places, round-precision=2]   
        S[table-text-alignment=right, scientific-notation=false, exponent-mode=fixed, round-mode=places, round-precision=2]   
        S[table-text-alignment=right, scientific-notation=false, exponent-mode=fixed, round-mode=places, round-precision=2]   
        |
    }
\hline
{$L_{\text{max}}$} & {$R_{\text{eff}}$} & {DOFs} & {$T_{\text{build}}$ (sec)} & {$T_{\text{upwards}}$ (sec)} & {$T_{\text{solve}}$ (sec)} & {$S$ (MB)} \\
\hline
\num{3} & \num{128} & \num{16384} & \num{3.01504249999999E-01} & \num{7.933354200E-02} & \num{4.67959999999999E-02} & \num{1.58822374343872E+01} \\
\num{4} & \num{256} & \num{65536} & \num{1.49072670900000E+00} & \num{3.993616660E-01} & \num{1.95166584000000E-01} & \num{8.15484972000122E+01} \\
\num{5} & \num{512} & \num{262144} & \num{7.21594766700000E+00} & \num{1.882675541E+00} & \num{9.14324582999999E-01} & \num{3.98233067512512E+02} \\
\num{6} & \num{1024} & \num{1048576} & \num{3.38984894589999E+01} & \num{8.795803542E+00} & \num{3.89846425000000E+00} & \num{1.88101041126251E+03} \\
\num{7} & \num{2048} & \num{4194304} & \num{1.59021219833999E+02} & \num{4.257397542E+01} & \num{1.88671120829999E+01} & \num{8.67619791126251E+03} \\
\hline
\num{3} & \num{128} & \num{8704} & \num{1.29889625000000E-01} & \num{2.823758300E-02} & \num{1.49769169999999E-02} & \num{4.63061237335205E+00} \\
\num{4} & \num{256} & \num{22528} & \num{3.62822042000000E-01} & \num{8.477891600E-02} & \num{2.33342910000000E-02} & \num{1.33817796707153E+01} \\
\num{5} & \num{512} & \num{54784} & \num{8.70019541999999E-01} & \num{1.982561250E-01} & \num{4.41358329999999E-02} & \num{2.84357728958129E+01} \\
\num{6} & \num{1024} & \num{163072} & \num{2.88925454100000E+00} & \num{7.241970830E-01} & \num{1.29217458000000E-01} & \num{1.16809662818908E+02} \\
\num{7} & \num{2048} & \num{485632} & \num{9.30737075000000E+00} & \num{2.389094875E+00} & \num{3.30502583999999E-01} & \num{3.95934556007385E+02} \\
\hline
    \end{tabular}
    \label{tab:helmholtz_timing}
\end{table}

%% file: sections/conclusion.tex
\section{Conclusion}
\label{sec:conclusion}

We have demonstrated a new implementation of the Hierarchical Poincaré-Steklov method on a quadtree-adaptive mesh. We outlined the key differences between a binary tree and quadtree structure, as well as the linear algebra associated with a quadtree implementation. Our novel full quadtree indexing allows for efficient storage of the data matrices needed in the HPS method.

Our numerical experiments show that the quadtree-adaptive HPS method is fast and efficient in terms of time, error, and data storage. We solved three elliptic partial differential equations to demonstrate the correctness and efficiency of our method. This kind of implementation is well-suited for an adaptive mesh framework such as \texttt{p4est}. Indeed, we have shown that the quadtree-adaptive HPS method can solve the elliptic problems to similar error with up to $17$ times speedup in the build stage and up to $57$ times speedup in the solve stage. The speedup depends on how well the mesh is able to adapt to the curvature of the solution. Memory use is of particular interest when solving elliptic PDEs with a direct method, and we have demonstrated that this is a highly efficient implementation.

A major advantage of this direct solver for elliptic problems is the ability to pre-compute the set of solution operators needed to solve the problem. Once the set of solution operators is known, a solve step can be done in linear time. This is especially powerful for time-stepping problems where one can pre-compute the solution operators, then apply them at subsequent time steps. Further development of this implementation is in progress to couple this elliptic solver to a hyperbolic finite volume solver through an operator splitting approach. In addition, a fully parallel implementation is being developed for use on large supercomputers to solve coupled hyperbolic/elliptic partial differential equations.